# Environmental management and restoration under unified risk and uncertainty using robustified dynamic Orlicz risk

Hidekazu Yoshioka[1, *], Motoh Tsujimura[2], Futoshi Aranishi[3] and Tomomi Tanaka[4]

[1] Japan Advanced Institute of Science and Technology, 1-1 Asahidai, Nomi, Ishikawa 923-1292, Japan

[2] Doshisha University, Karasuma-Higashi-iru, Imadegawa-dori, Kamigyo-ku, Kyoto 602-8580, Japan

[3] Shimane University Institute of Agricultural and Life Sciences, 1060 Nishikawatsu, Matsue, Shimane 690-8504, Japan

[4] Shimane University, 1060 Nishikawatsu, Matsue, Shimane 690-8504, Japan

[*] Corresponding author: yoshih@jaist.ac.jp

ORCID:0000-0002-5293-3246 (Hidekazu Yoshioka), 0000-0001-6975-9304 (Motoh Tsujimura), 0000-0002-6463-3275 (Futoshi Aranishi)

**Abstract** Environmental management and restoration should be designed such that the risk and uncertainty owing to nonlinear stochastic systems can be successfully addressed. We apply the robustified dynamic Orlicz risk to the modeling and analysis of environmental management and restoration to consider both the risk and uncertainty within a unified theory. We focus on the control of a jump-driven hybrid stochastic system that represents macrophyte dynamics. The dynamic programming equation based on the Orlicz risk is first obtained heuristically, from which the associated Hamilton–Jacobi–Bellman (HJB) equation is derived. In the proposed Orlicz risk, the risk aversion of the decision-maker is represented by a power coefficient that resembles a certainty equivalence, whereas the uncertainty aversion is represented by the Kullback–Leibler divergence, in which the risk and uncertainty are handled consistently and separately. The HJB equation includes a new state-dependent discount factor that arises from the uncertainty aversion, which leads to a unique, nonlinear, and nonlocal term. The link between the proposed and classical stochastic control problems is discussed with a focus on control-dependent discount rates. We propose a finite difference method for computing the HJB equation. Finally, the proposed model is applied to an optimal harvesting problem for macrophytes in a brackish lake that contains both growing and drifting populations.

## 1. Introduction
### 1.1 Problem background

Environmental management and restoration are necessary for the sustainable coexistence of nature and humans. Water environments constitute a principal part of aquatic ecosystems and socioenvironmental systems [1–2]. Therefore, they are important research targets in civil and environmental engineering, and related fields [3–5]. Environmental management and restoration are decision-making problems for complex systems because they deal with highly unpredictable dynamics that involve physical, biological, chemical, and social phenomena [6–8]. Mathematical models can provide useful insights into the design of cost-effective environmental management and restoration schemes by decision-makers.

Environmental dynamics can be described using stochastic differential equations (SDEs), which are (ordinary) differential equations that are perturbed by random processes. SDEs and related mathematical models are versatile tools for describing various phenomena, including renewable resource dynamics [9], the transport of water [10] and radionuclides [11], algal patch dynamics [12], and the temporal evolution of the water temperature [13]. Stochastic control [14] is theoretically and practically effective in the optimization and control of environmental dynamics that are modelled by SDEs. The main task of stochastic control is to solve the optimality equation, which is often known as the Hamilton–Jacobi–Bellman (HJB) equation, based on an objective function to be optimized subject to the system dynamics to be controlled. Applied problems in this area include waste cleaning by a mining firm [15], the harvesting of fishery resources in marine environments [16], the planning of mitigation measures for biodiversity and connectivity losses owing to traffic networks [17], lake pollution management [18], biofilm growth [19], and the control of dam-reservoir systems [20].

Environmental management and restoration face risk and uncertainty owing to stochastic system dynamics. The risk here refers to extreme events that may occur with low probabilities, such as cyanobacterial blooms [21] and severe drought [22]. The uncertainty arises from incomplete knowledge regarding the target system to be managed owing to poor and/or sparse data availability, which emerges as ambiguous parameter values [23–24] and functional shapes of coefficients [25–26] in the system dynamics. The control of the risk and uncertainty in planning environmental management and restoration schemes has been addressed from the perspective of static decision theories [27–28] as well as dynamic theories based on stochastic control under a distorted measure [29–30] and nonlinear certainty equivalence [31–32]. However, to the best of our knowledge, a unified theory for dealing with the risk and uncertainty concurrently remains lacking.

The modern theory for decision-making subject to risk and uncertainty in mathematical insurance and finance suggests that robustified risk measures that are defined through Orlicz norms may serve as the objective function in the optimization problem [33–35]. The Orlicz norm is a generation of the Lebesgue norm that evaluates risks more flexibly by covering a wider range of functional variability [36]. These risk measures are well defined in certain Orlicz spaces that are Banach spaces and can formulate the risk and uncertainty in a consistent and separable manner. The risk that represents extremely small or large values of the random variable of interest is evaluated through the Young function (sometimes known as the Orlicz function), which is understood as the risk aversion coefficient, in the Orlicz norm. The uncertainty is accounted for by a penalty function such as the divergence (e.g., [37–38]). Time-consistent dynamic Orlicz premia that resemble dynamic programming in stochastic control were recently proposed, and their mathematical properties were studied in detail [39], suggesting their versatility in applications; however, to the best of our knowledge, such approaches are still rare.

**1.2 Objectives and contributions**

The objectives of this study are to propose and analyze a robustified dynamic Orlicz risk, which is known as the robust Orlicz risk, for the stochastic control of nonlinear and hybrid SDEs that arise in environmental management and restoration over an infinite horizon. The SDEs are jump-driven and represent the population dynamics of aquatic species such as macrophytes that are subject to time-discrete intervention by an environmental manager, namely the decision-maker (DM), and catastrophic events such as floods. These assumptions simplify the target problem but still highlight the role of the Orlicz risk in stochastic control subject to the risk and uncertainty. The problem in this study is related to many real problems, such as the sustainable conservation of biodiversity through macrophyte removal [40] and the prediction/control of algal blooms [41–42].

The robust Orlicz risk is defined using a certainty equivalent of the power function that represents the risk aversion and the Kullback–Leibler divergence [43] for measuring the uncertainty aversion of frequent use in optimization and learning research. The entropic formalism is advantageous from a practical perspective, as the HJB equation that is associated with the robust Orlicz risk admits analytical solutions in certain cases despite being highly nonlinear. Our approach provides optimal control considering the uncertainty.

The proposed dynamic programming approach is related to stochastic control based on the Epstein–Zin nonlinear expectation [44–47]. However, our model assumes a more general coefficient shape and accounts for the uncertainty. We demonstrate that the HJB equation possesses a unique nonlinear and nonlocal term that arises from a state- and control-dependent discount rate. Although state-dependent discount rates [48–49] and state- and control-dependent discount rates [50] have been discussed in the context of Markov decision processes, they have not been integrated into the robust Orlicz risk. The use of the HJB equation enables us to obtain the optimal controls and further evaluate their performance.

In addition to mathematically analyzing the HJB equation, we propose a simple finite difference method for its discretization using common one-sided and mid-point differences (e.g., [51–53]). The monotonicity and stability properties of the discretization are crucial for obtaining convergent numerical solutions to the HJB and related differential equations [54–57]. Finally, we apply the proposed Orlicz risk and HJB equation to the harvesting problem of drifting aquatic vegetation in a shallow lake.

The remainder of this paper is organized as follows. Section 2 introduces the mathematical model, including the SDE to be controlled, the robust Orlicz risk, and the HJB equation. In Section 3, we present the mathematical analysis results of the model and investigated the finite difference method for the HJB equation. Section 4 describes the application of the HJB equation. The conclusions and future perspectives are provided in Section 5. Appendix A presents the proofs of Propositions. Appendices B and C provide auxiliary results that are not presented in the main text.

## 2. Mathematical model
### 2.1 Target system in absence of uncertainty

We focus on a two-population system of macrophyte dynamics driven by jump noise processes that represent environmental perturbations and human interventions. The system leads to an exactly solvable case under certain assumptions, and it contains nontrivial features that are subject to risk and uncertainty.

The stochastic system dynamics are defined in a complete probability space $(\Omega, \mathcal{F}, \mathbb{P})$ as usual, and are expressed as follows in the absence of uncertainty:

$$\underbrace{\mathrm{d}\begin{pmatrix} X_t^{(1)} \\ X_t^{(2)} \end{pmatrix}}_{\text{Increment}} = \underbrace{\begin{pmatrix} a_1\left(X_t^{(1)}, X_t^{(2)}\right) \\ a_2\left(X_t^{(1)}, X_t^{(2)}\right) \end{pmatrix} \mathrm{d}t}_{\text{Deterministic growth/decay}} + \underbrace{\begin{pmatrix} b_1\left(X_{t-}^{(1)}, X_{t-}^{(2)}\right) \\ b_2\left(X_{t-}^{(1)}, X_{t-}^{(2)}\right) \end{pmatrix} \mathrm{d}N_t}_{\text{Stochastic growth/decay}} + \underbrace{\begin{pmatrix} c_1\left(X_{t-}^{(1)}, X_{t-}^{(2)}, u_t^{(1)}\right) \\ c_2\left(X_{t-}^{(1)}, X_{t-}^{(2)}, u_t^{(2)}\right) \end{pmatrix} \mathrm{d}Z_t}_{\text{Human intervention}}, \quad t > 0, \qquad (1)$$

where $t \geq 0$ is the time, $t-$ represents the left limit at $t$, $X^{(i)} = \left(X_t^{(i)}\right)_{t \geq 0}$ is the $i$ th population as a nonnegative continuous-time stochastic process ($i = 1, 2$), $a_i : [0, +\infty)^2 \to \mathbb{R}$ is the drift coefficient for the $i$ th population ($i = 1, 2$), $b_i : [0, +\infty)^2 \to \mathbb{R}$ is the uncontrollable perturbation of the $i$ th population ($i = 1, 2$) at exogenous jumps that are generated by a counting process $N = (N_t)_{t \geq 0}$, $c_i : [0, +\infty)^2 \times U \to C$ with state-independent compact sets $U, C$, to simplify the problem, is the controllable increase or decrease in the $i$ th population ($i = 1, 2$) by another counting process $Z = (Z_t)_{t \geq 0}$ with an increment of 1 at each jump, and $u^{(i)} = \left(u_t^{(i)}\right)_{t \geq 0}$ is the control process by the DM. Only the values of $u^{(i)}$ at the jumps of $Z$ affect the dynamics (1). The counting processes $N$ and $Z$ are assumed to be mutually independent Poisson processes with constant intensities $\lambda_N > 0$ and $\lambda_Z > 0$, respectively. The jump size of $N$ follows a probability density function $P$ with a compact support in $[-1, +\infty)$. The support is denoted as $\mathrm{Supp}P = [\underline{z}, \overline{z}]$ unless otherwise specified, with constants $\underline{z}, \overline{z}$ satisfying $-1 < \underline{z} < \overline{z} < +\infty$. Because a fixed scheduling planned by the DM may be accidentally randomized due to unexpectable reasons such as a climatic condition, human interventions are assumed to occur at random times.

At each jump time $s$ of $N$ with jump size $\Delta N_s$, the populations are perturbed as

$$X_s^{(i)} = X_{s-}^{(i)} + b_i\left(X_{s-}^{(1)}, X_{s-}^{(2)}\right) \Delta N_s \quad (i = 1, 2). \qquad (2)$$

Similarly, at each jump time $s$ of $Z$, the populations are perturbed as

$$X_s^{(i)} = X_{s-}^{(i)} + c_i\left(X_{s-}^{(1)}, X_{s-}^{(2)}, u_s^{(i)}\right) \quad (i = 1, 2). \qquad (3)$$

We assume that the coefficients $a_i, b_i, c_i$ are sufficiently regular so that the system (1) starting from a nonnegative and deterministic initial condition $\left(X_0^{(1)}, X_0^{(2)}\right)$ admits a unique nonnegative càdlàg solution $\left(X_t^{(1)}, X_t^{(2)}\right)_{t \geq 0}$.

The following assumption is used throughout the paper unless otherwise specified.

***Assumption 1***

1. Both $a_1, a_2$ are globally Lipschitz continuous in $[0, +\infty)^2$, and there exists a constant $\theta > 0$ such that $a_1(x_1, x_2), a_2(x_1, x_2) < 0$ for $x_1, x_2 \geq \theta$. Moreover, $a_1(0, \cdot), a_2(\cdot, 0) \geq 0$.

2. Both $b_1, b_2$ are globally Lipschitz continuous and uniformly bounded in $[0, +\infty)^2$. Moreover, in $[0, +\infty)^2$ they satisfy $x_i + \underline{z} b_i(x_1, x_2) \geq 0$ as well as $x_i + \overline{z} b_i(x_1, x_2) \leq \overline{b}$ ($i = 1, 2$) with a positive constant $\overline{b} > 0$.

3. Both $c_1, c_2$ are globally Lipschitz continuous and uniformly bounded in $[0, +\infty)^2 \times C$. Moreover, in $[0, +\infty)^2 \times C$ they satisfy $0 \leq x_i + c_i(x_1, x_2, u_i) \leq \overline{c}$ ($i = 1, 2$) with a positive constant $\overline{c} > 0$.

Assumption 1 guarantees that the state variables are constrained in the compact set $\left[0, \max\{\theta, \overline{b}, \overline{c}\}\right]^2$ almost surely if the initial state $\left(X_0^{(1)}, X_0^{(2)}\right)$ is sufficiently close to $(0, 0)$. Furthermore, Assumption 1.1 shows that SDE (1) in the absence of jumps admits a unique pathwise solution that is bounded globally in time. Assumptions 1.2-1.3 demonstrate that the uniqueness is still preserved under jump disturbance. Assumption 1 is satisfied in our application study presented later.

***Remark 1*** The harvesting problem of aquatic vegetation consisting of both growing ($X^{(1)}$) and drifting ($X^{(2)}$) populations is an example of coefficients that comply with Assumption 1 (in some compact set). This is a two-state case in which the dynamics of the two population patches, namely the source and sink patches, are coupled, and $t$ is a modification of the case considered in [58]. The populations should be maintained at moderate levels from ecological and fisheries perspectives:

$$a_i(x_1, x_2) = \chi_\mu(x_i)\left\{r_i x_i\left(1 - \frac{x_i}{Q_i}\right) - d_i x_i + e_i x_{3-i}\right\} \quad \text{(logistic-like growth with immigration)} \quad (4)$$

for small $x_1, x_2$, with constants $r_i, Q_i > 0$ and $d_i, e_i \geq 0$ in which $r_i > d_i$, and a truncation $\mu > 0$:

$$\chi_\mu(x_i) = \min\left\{1, \frac{1 - x_i / Q_i}{\mu}\right\}, \quad (5)$$

$$b_i(x_1, x_2) = -x_i \quad \text{(proportional decay)}, \quad (6)$$

$$c_i(x_1, x_2, u_i) = -\min\{x_i, u_i\} \quad \text{(linear harvesting with saturation)}. \quad (7)$$

The truncation function (5) is considered such that the population is within the carrying capacity $Q_i$. The influence of the truncation decreases with a smaller $\mu > 0$.

***Remark 2*** We can consider a single-population system and a higher-dimensional system by suitably modifying (1). Moreover, diffusive noise terms that are driven by Brownian motion can be added to the dynamics. For theoretical simplicity, we only consider the jump-driven dynamics. One may add a Brownian noise to the dynamics, which would add a second-order partial differential term in the HJB equation.

## 2.2 System subject to uncertainty

We consider that the intensity and size of the jumps of process $N$ are uncertain for the DM. The uncertainty is the perturbation of the jump size and intensity in terms of the Radon–Nikodým derivative [59–60]. The uncertainty is represented by a process that is adapted to a natural filtration $\mathcal{F}$ generated by the two jump processes $N, Z$. We set the compensated process $\tilde{N}_t := N_t - \lambda_N \int_0^t \left( \int_{-1}^{+\infty} P(z) \mathrm{d}z \right) \mathrm{d}t$ ($t \geq 0$), which is a martingale under the benchmark probability measure $\mathbb{P}$: the probability measure in the absence of uncertainty.

Furthermore, we consider a distorted probability measure $\mathbb{P}(\phi)$ that is equivalent to $\mathbb{P}$ under which the compensated process $\tilde{N}_t^\phi := N_t - \lambda_N \int_0^t \left( \int_{-1}^{+\infty} \phi_t(z) P(z) \mathrm{d}z \right) \mathrm{d}t$ becomes a martingale, where $\phi = \left( \phi_t(\cdot) \right)_{t \geq 0}$ is a positive random field adapted to $\mathcal{F}$. We obtain $\tilde{N}_t^\phi = \tilde{N}_t$ if $\phi_t(\cdot) \equiv 1$, which means that no uncertainty exists; that is, $\mathbb{P} = \mathbb{P}(\phi)\big|_{\phi \equiv 1}$. The expectation under $\mathbb{P}(\phi)$ is denoted as $\mathbb{E}_{\mathbb{P}(\phi)}$. We consider the SDE (1) under a distorted probability measure for each $\phi$ that belongs to some admissible set. The process $\phi$ is viewed as a control variable to be optimized by an opponent of the DM (so-called nature). Thus, we formulate the control problem as a zero-sum differential game.

***Remark 3*** The other counting process $Z$ can also be made uncertain. The proposed Orlicz risk converts the uncertainty in process $N$ into a nonlinear and nonlocal term in the HJB equation. This correspondence may be the same for $Z$ if its intensity and size are uncertain. Therefore, we only consider the uncertainty of $N$ in the theoretical analysis.

## 2.3 Robust Orlicz risk
### 2.3.1 Abstract formulation

We define the Orlicz risk [33] as a type of recursive utility. We first briefly review Orlicz spaces and norms (e.g., [36]). The Orlicz space is a generalization of Lebesgue spaces, which are often known as $L^p$ spaces. In the following, $\mathbb{E}$ represents the generic expectations and $Y$ is a nonnegative random variable that is adapted to $\mathcal{F}$. Given an Orlicz function $\Phi:[0,+\infty)\to[0,+\infty)$ that is continuous, increasing, and convex with $\Phi(0)=0$ and $\Phi(1)=1$, the Orlicz space $L_\Phi$ is defined as

$$L_\Phi := \left\{Y : \mathbb{E}\left[\Phi\left(\frac{Y}{h}\right)\right] < +\infty \text{ for some } h > 0\right\} \tag{8}$$

with the Orlicz norm

$$\|Y\|_\Phi := \inf\left\{h > 0 : \mathbb{E}\left[\Phi\left(\frac{Y}{h}\right)\right] \leq 1\right\}. \tag{9}$$

If $\Phi(x)=x^p$ ($p\geq 1$), the Orlicz space is identified as the Lebesgue space by the norm equivalence

$$\|Y\|_p := \left(\mathbb{E}[Y^p]\right)^{\frac{1}{p}} = \|Y\|_\Phi. \tag{10}$$

The dynamic Orlicz risk was defined in [39] as

$$\|Y_t\|_\Phi := \inf\left\{h_t > 0, h \text{ is bounded and adapted to } \mathcal{F} : \mathbb{E}\left[\Phi\left(\frac{Y_t}{h_t}\right)\bigg|\mathcal{F}_t\right] \leq 1\right\}, \quad t \geq 0. \tag{11}$$

In this context, $\mathbb{E}=\mathbb{E}_\mathbb{P}$ in (11). A robustified version was also proposed [39]; for $t \geq 0$, set

$$\|Y_t\|_{c,\Phi} := \inf\left\{h_t > 0, h \text{ is bounded and adapted to } \mathcal{F} : \sup_\phi\left(\mathbb{E}_{\mathbb{P}(\phi)}\left[\Phi\left(\frac{Y_t}{h_t}\right)\bigg|\mathcal{F}_t\right] - c_t(\phi)\right) \leq 1\right\}, \tag{12}$$

with a nonnegative functional $c_t(\phi)$ of the positive and measurable process $\phi$ up to time $t$, such that $c_t(\phi)|_{\phi\equiv 1} = 0$. If $\Phi(x)=x^p$ ($p\geq 1$), (12) becomes

$$\|Y_t\|_{c,\Phi} = \sup_\phi\left\{\frac{\mathbb{E}_{\mathbb{P}(\phi)}\left[(Y_t)^p\big|\mathcal{F}_t\right]}{1+c_t(\phi)}\right\}^{\frac{1}{p}}, \quad t \geq 0. \tag{13}$$

This formulation resembles the certainty equivalent (10) but with a discount $(1+c_t(\phi))^{-1} \leq 1$. Thus, the dynamic Orlicz risk (13) controls the risk by selecting the power index $p \geq 1$ so that a larger risk aversion (i.e., sensitivity against larger objective function value) of the DM is assumed for a larger $p$. The uncertainty may distort the expectation $\mathbb{E}_{\mathbb{P}(\phi)}$ through the process $\phi$, and its impacts can be bounded by the existence of the discount factor $(1+c_t(\phi))^{-1}$. To the best of our knowledge, this type of formulation has not yet been discussed. The right-hand side of (13) is time consistent if $\mathbb{E}_{\mathbb{P}(\phi)}\left[\Phi\left(\frac{Y_t}{h_t}\right)\bigg|\mathcal{F}_t\right] - c_t(\phi)$ is (Lemma 18 of [39]).

## 2.4 Recursive formulation

We heuristically obtain a recursive formulation of a stochastic control problem. We let $k > 0$ be a small time increment that is later set to 0. We fix the time as $t > 0$. We consider that both costs and disutility arise in the time interval $(t, t+k)$, provided that the information at $t$ is given. The formulation is heuristic and qualitatively follows the literature [61–62].

There is a cost of intervention by the DM, which may occur at each jump of $Z$. We assume that the instantaneous cost incurred by $(u^{(1)}, u^{(2)}) \in C \times C$ is $g(u^{(1)}, u^{(2)})$ with a nonnegative bounded function $g : C \times C \to [0, +\infty)$ such that $g(0, 0) = 0$. The total cost incurred during $(t, t+k)$ is

$$\sum_{s \in (t, t+k), Z} g\left(u_s^{(1)}, u_s^{(2)}\right), \tag{14}$$

where the summation $\sum_{s \in (t_1, t_2), Z}$ contains the summands at all jumps of $Z$ in $(t_1, t_2)$. The summation for the process $N$ is defined similarly. We also consider a disutility caused by the deviation of the population states from prescribed target values. Subsequently, using a nonnegative bounded and uniformly continuous function $f : [0, +\infty)^2 \to [0, +\infty)$, the disutility during $(t, t+k)$ is set to

$$\int_t^{t+k} f\left(X_s^{(1)}, X_s^{(2)}\right) \mathrm{d}s. \tag{15}$$

Furthermore, we need to specify the penalization $c_t$ against uncertainty. As the process $\phi$ is viewed as a Radon–Nikodým derivative between the two probability measures $\mathbb{P}$ and $\mathbb{P}(\phi)$, we use the divergence to quantify their difference in terms of $\phi$ (e.g., [37–38]):

$$c_t(\phi) = \frac{\lambda_N}{\psi} \mathbb{E}_{\mathbb{P}(\phi)} \left[ \int_t^{t+k} \left( \int_{-1}^{+\infty} D(\phi_s(z)) P(z) \mathrm{d}z \right) \mathrm{d}s \,\middle|\, \mathcal{F}_t \right], \quad D(\phi_s) := \phi_s \ln \phi_s - \phi_s + 1, \tag{16}$$

where $\psi > 0$ is an uncertainty aversion parameter such that the DM is more uncertainty averse with a larger $\psi$. The right-hand side of (16) is understood as the mean cumulative relative entropy during $(t, t+k)$. Moreover, $D$ is nonnegative, convex, and has the global minimum $D(1) = 0$. Thus, a positive divergence value means that there is some uncertainty in the model.

The admissible sets of controls must also be specified. The admissible set $\mathcal{A}$ of $(u^{(1)}, u^{(2)})$ is determined as follows (e.g., [58,63]):

$$\mathcal{A} = \left\{ (u^{(1)}, u^{(2)}) \,\middle|\, \begin{array}{l} (u_t^{(1)}, u_t^{(2)}) \in C \times C \text{ and is measurable for } t > 0, \\ \text{and } (u_t^{(1)}, u_t^{(2)}) = (u_{\tau_{i-1}}^{(1)}, u_{\tau_{i-1}}^{(2)}) \ (\tau_{i-1} \leq t < \tau_i) \text{ for } i \in \mathbb{N} \end{array} \right\}, \tag{17}$$

where $\tau_i$ ($i = 0, 1, 2, ...$) is the $i$th jump time of process $Z$ with $\tau_0 = 0$ for brevity, which constitutes an increasing sequence. We obtain $\tau_{i-1} < \tau_i$ ($i = 1, 2, 3...$) almost surely owing to the assumption that $Z$ is a compound Poisson process, the jump times of which differ from those of $N$ with probability 1. Only the values of $(u^{(1)}, u^{(2)})$ at each $\tau_i$ are important as implied in (3).

The admissible set $\mathcal{B}$ of $\phi$ is determined as follows (e.g., [64]):

$$\mathcal{B} = \{\phi | \phi_t(\cdot) > 0 \text{ and is predictable for } t > 0, \text{ and satisfies (19)-(20)}\}, \tag{18}$$

$$\int_0^t \left( \int_{-1}^{+\infty} D(\phi_s(z)) P(z) dz \right) ds < +\infty \quad \text{for } t > 0, \tag{19}$$

$$\mathbb{E}_\mathbb{P}\left[ \exp\left( \frac{1}{2} \int_0^t \int_{-1}^{+\infty} \{\phi_s(z) \ln \phi_s(z) - \phi_s(z) + 1\} P(z) dz ds \right) \right] < +\infty \quad \text{for } t > 0. \tag{20}$$

Sets $\mathcal{A}, \mathcal{B}$ are nonempty because they contain the null controls $(u^{(1)}, u^{(2)}) \equiv (0,0)$ and $\phi \equiv 1$, respectively, which are not indicated in the equations below to avoid confusion.

We set a constant discount rate $\delta > 0$, as in classical control problems, over an infinite horizon. Based on the Orlicz risk, we set the recursive utility process $\Psi = (\Psi_t)_{t \geq 0}$ as follows:

$$\Psi_t = \inf_{\mathbf{u}} \sup_\phi \left\{ \overbrace{\mathbb{E}_{\mathbb{P}(\phi)}\left[ \int_t^{t+k} f^p(X_s^{(1)}, X_s^{(2)}) ds + \sum_{s \in (t, t+k), Z} g^p(u_s^{(1)}, u_s^{(2)}) \bigg| \mathcal{F}_t \right]}^{\text{Cumulative disutility and cost}} + e^{-\delta k} \underbrace{\frac{\mathbb{E}_{\mathbb{P}(\phi)}[(\Psi_{t+k})^p | \mathcal{F}_t]}{1 + \lambda_N \psi^{-1} \mathbb{E}_{\mathbb{P}(\phi)}\left[ \int_t^{t+k} \left( \int_{-1}^{+\infty} D(\phi_s(z)) P(z) dz \right) ds \bigg| \mathcal{F}_t \right]}}_{\text{Discounted future value}} \right\}^{\frac{1}{p}}, \quad t > 0. \tag{21}$$

(Current value on the left-hand side)

Heuristically, equation (21) recursively determines $\Psi$ backward over time. Here, $\Psi$ represents a process that accumulates the cost and disutility under the aversion against the risk and uncertainty, and the recursion (21) computes the current value as a weighted mean of the cumulative disutility and cost and the discounted future value, as in classical dynamic programming, but in a nonlinear and hence risk-averse manner. The discount rate $\delta$ must be sufficiently high in some cases to define $\Psi$ effectively. Finally, taking the $p$th power of $f, g$ guarantees the consistency of units among the utility, cost, and disutility.

*Remark 4* Some utility by the existence of the population can also be considered in the future, such as the bioeconomic use of harvested vegetation [65].

## 2.5 HJB equation

As Markovian control is natural in both theory and application, we seek $\Psi$ leading to Markovian controls. With an abuse of notation, we write $\Psi_t = \Psi(X_t^{(1)}, X_t^{(2)})$ and refer to $\Psi : [0, +\infty)^2 \to \mathbb{R}$ on the right-hand side as the value function as in classical control theory. If $(X_t^{(1)}, X_t^{(2)}) = (x_1, x_2)$ at $t$, (21) yields

$$\left(\Psi\left(x_{1},x_{2}\right)\right)^{p}=\inf_{\mathbf{u}}\sup_{\phi}\left\{\begin{array}{l}\mathbb{E}_{\mathbb{P}(\phi)}\left[\int_{t}^{t+k}f^{p}\left(X_{s}^{(1)},X_{s}^{(2)}\right)\mathrm{d}s+\sum_{s\in(t,t+k),Z}g^{p}\left(u_{s}^{(1)},u_{s}^{(2)}\right)\bigg|\mathcal{F}_{t}\right]\\ +e^{-\delta k}\dfrac{\mathbb{E}_{\mathbb{P}(\phi)}\left[\left(\Psi\left(X_{t+k}^{(1)},X_{t+k}^{(2)}\right)\right)^{p}\bigg|\mathcal{F}_{t}\right]}{1+\lambda_{N}\psi^{-1}\mathbb{E}_{\mathbb{P}(\phi)}\left[\int_{t}^{t+k}\left(\int_{-1}^{+\infty}D(\phi_{s}(z))P(z)\mathrm{d}z\right)\mathrm{d}s\bigg|\mathcal{F}_{t}\right]}\end{array}\right\},\ t>0. \quad (22)$$

We will use $F=\Psi^{p}$ to simplify the notation. From (22), we obtain

$$0=\frac{1}{k}\inf_{\mathbf{u}}\sup_{\phi}\left\{\begin{array}{l}\mathbb{E}_{\mathbb{P}(\phi)}\left[\int_{t}^{t+k}f^{p}\left(X_{s}^{(1)},X_{s}^{(2)}\right)\mathrm{d}s+\sum_{s\in(t,t+k),Z}g^{p}\left(u_{s}^{(1)},u_{s}^{(2)}\right)\bigg|\mathcal{F}_{t}\right]\\ +e^{-\delta k}\dfrac{\mathbb{E}_{\mathbb{P}(\phi)}\left[\left(\Psi\left(X_{t+k}^{(1)},X_{t+k}^{(2)}\right)\right)^{p}\bigg|\mathcal{F}_{t}\right]}{1+\lambda_{N}\psi^{-1}\mathbb{E}_{\mathbb{P}(\phi)}\left[\int_{t}^{t+k}\left(\int_{-1}^{+\infty}D(\phi_{s}(z))P(z)\mathrm{d}z\right)\mathrm{d}s\bigg|\mathcal{F}_{t}\right]}-\left(\Psi\left(x_{1},x_{2}\right)\right)^{p}\end{array}\right\}$$

$$=\inf_{\mathbf{u}}\sup_{\phi}\left\{\begin{array}{l}\dfrac{1}{k}\mathbb{E}_{\mathbb{P}(\phi)}\left[\int_{t}^{t+k}f^{p}\left(X_{s}^{(1)},X_{s}^{(2)}\right)\mathrm{d}s+\sum_{s\in(t,t+k),Z}g^{p}\left(u_{s}^{(1)},u_{s}^{(2)}\right)\bigg|\mathcal{F}_{t}\right]\\ +\dfrac{1}{k}e^{-\delta k}\dfrac{\mathbb{E}_{\mathbb{P}(\phi)}\left[\left(\Psi\left(X_{t+k}^{(1)},X_{t+k}^{(2)}\right)\right)^{p}\bigg|\mathcal{F}_{t}\right]-e^{\delta k}\left(\Psi\left(x_{1},x_{2}\right)\right)^{p}}{1+\lambda_{N}\psi^{-1}\mathbb{E}_{\mathbb{P}(\phi)}\left[\int_{t}^{t+k}\left(\int_{-1}^{+\infty}D(\phi_{s}(z))P(z)\mathrm{d}z\right)\mathrm{d}s\bigg|\mathcal{F}_{t}\right]}\\ -\dfrac{1}{k}\dfrac{\lambda_{N}\psi^{-1}\mathbb{E}_{\mathbb{P}(\phi)}\left[\int_{t}^{t+k}\left(\int_{-1}^{+\infty}D(\phi_{s}(z))P(z)\mathrm{d}z\right)\mathrm{d}s\bigg|\mathcal{F}_{t}\right]}{1+\lambda_{N}\psi^{-1}\mathbb{E}_{\mathbb{P}(\phi)}\left[\int_{t}^{t+k}\left(\int_{-1}^{+\infty}D(\phi_{s}(z))P(z)\mathrm{d}z\right)\mathrm{d}s\bigg|\mathcal{F}_{t}\right]}\left(\Psi\left(x_{1},x_{2}\right)\right)^{p}\end{array}\right\}. \quad (23)$$

Furthermore, we obtain

$$\lim_{k\to+0}\frac{1}{k}\mathbb{E}_{\mathbb{P}(\phi)}\left[\int_{t}^{t+k}f^{p}\left(X_{s}^{(1)},X_{s}^{(2)}\right)\mathrm{d}s+\sum_{s\in(t,t+k),Z}g^{p}\left(u_{s}^{(1)},u_{s}^{(2)}\right)\bigg|\mathcal{F}_{t}\right]=f^{p}\left(x_{1},x_{2}\right)+\lambda_{Z}g^{p}\left(u_{t}^{(1)},u_{t}^{(2)}\right), \quad (24)$$

$$\lim_{k\to+0}\frac{1}{k}\mathbb{E}_{\mathbb{P}(\phi)}\left[\int_{t}^{t+k}\left(\int_{-1}^{+\infty}D(\phi_{s}(z))P(z)\mathrm{d}z\right)\mathrm{d}s\bigg|\mathcal{F}_{t}\right]=\int_{-1}^{+\infty}D(\phi_{t}(z))P(z)\mathrm{d}z, \quad (25)$$

and by Itô's formula for jump processes (e.g., Theorem 1.14 of [14]),

$$\lim_{k\to+0}\frac{\mathbb{E}_{\mathbb{P}(\phi)}\left[\left(\Psi\left(X_{t+k}^{(1)},X_{t+k}^{(2)}\right)\right)^{p}\Big|\mathcal{F}_{t}\right]-e^{\delta k}\left(\Psi(x_{1},x_{2})\right)^{p}}{k}$$

$$=\lim_{k\to+0}\frac{\mathbb{E}_{\mathbb{P}(\phi)}\left[\left(\Psi\left(X_{t+k}^{(1)},X_{t+k}^{(2)}\right)\right)^{p}\Big|\mathcal{F}_{t}\right]-\left(\Psi(x_{1},x_{2})\right)^{p}-\left(e^{\delta k}-1\right)\left(\Psi(x_{1},x_{2})\right)^{p}}{k}$$

$$=-\delta F(x_{1},x_{2})+a_{1}(x_{1},x_{2})\frac{\partial F(x_{1},x_{2})}{\partial x_{1}}+a_{2}(x_{1},x_{2})\frac{\partial F(x_{1},x_{2})}{\partial x_{2}}$$
$$+\lambda_{Z}\left(F\left(x_{1}+c_{1}(x_{1},x_{2},u_{t}^{(1)}),x_{2}+c_{2}(x_{1},x_{2},u_{t}^{(1)})\right)-F(x_{1},x_{2})\right)$$
$$+\lambda_{N}\int_{-1}^{+\infty}\left(F(x_{1}+b_{1}(x_{1},x_{2})z,x_{2}+b_{2}(x_{1},x_{2})z)-F(x_{1},x_{2})\right)\phi_{t}(z)P(z)\mathrm{d}z$$
(26)

Consequently, we obtain

$$0=\lim_{k\to+0}\frac{1}{k}\inf_{\mathbf{u}}\sup_{\phi}\left\{\begin{array}{l}\mathbb{E}_{\mathbb{P}(\phi)}\left[\int_{t}^{t+k}f^{p}\left(X_{s}^{(1)},X_{s}^{(2)}\right)\mathrm{d}s+\sum_{s\in(t,t+k),Z}g^{p}\left(u_{s}^{(1)},u_{s}^{(2)}\right)\Big|\mathcal{F}_{t}\right]\\ +e^{-\delta k}\dfrac{\mathbb{E}_{\mathbb{P}(\phi)}\left[\left(\Psi\left(X_{t+k}^{(1)},X_{t+k}^{(2)}\right)\right)^{p}\Big|\mathcal{F}_{t}\right]}{1+\lambda_{N}\psi^{-1}\mathbb{E}_{\mathbb{P}(\phi)}\left[\int_{t}^{t+k}\left(\int_{-1}^{+\infty}D(\phi_{s}(z))P(z)\mathrm{d}z\right)\mathrm{d}s\Big|\mathcal{F}_{t}\right]}-\left(\Psi(x_{1},x_{2})\right)^{p}\end{array}\right\}$$

$$=\inf_{\left(u_{t}^{(1)},u_{t}^{(2)}\right)}\sup_{\phi_{t}(\cdot)>0}\left\{\begin{array}{l}-\left(\delta+\lambda_{N}\psi^{-1}\int_{-1}^{+\infty}D(\phi(z))p(z)\mathrm{d}z\right)F(x_{1},x_{2})+f^{p}(x_{1},x_{2})+\lambda_{Z}g^{p}\left(u_{t}^{(1)},u_{t}^{(2)}\right)\\ +a_{1}(x_{1},x_{2})\dfrac{\partial F(x_{1},x_{2})}{\partial x_{1}}+a_{2}(x_{1},x_{2})\dfrac{\partial F(x_{1},x_{2})}{\partial x_{2}}\\ +\lambda_{Z}\left(F\left(x_{1}+c_{1}(x_{1},x_{2},u_{t}^{(1)}),x_{2}+c_{2}(x_{1},x_{2},u_{t}^{(1)})\right)-F(x_{1},x_{2})\right)\\ +\lambda_{N}\int_{-1}^{+\infty}\left(F(x_{1}+b_{1}(x_{1},x_{2})z,x_{2}+b_{2}(x_{1},x_{2})z)-F(x_{1},x_{2})\right)\phi_{t}(z)P(z)\mathrm{d}z\end{array}\right\}.$$
(27)

That is, we arrive at the HJB equation

$$\inf_{(u_{1},u_{2})\in C\times C}\sup_{\phi(\cdot)>0}\left\{\begin{array}{l}f^{p}(x_{1},x_{2})-\left(\delta+\lambda_{N}\psi^{-1}\int_{-1}^{+\infty}D(\phi(z))P(z)\mathrm{d}z\right)F(x_{1},x_{2})\\ +a_{1}(x_{1},x_{2})\dfrac{\partial F(x_{1},x_{2})}{\partial x_{1}}+a_{2}(x_{1},x_{2})\dfrac{\partial F(x_{1},x_{2})}{\partial x_{2}}\\ +\lambda_{Z}\left(F(x_{1}+c_{1}(x_{1},x_{2},u_{1}),x_{2}+c_{2}(x_{1},x_{2},u_{2}))-F(x_{1},x_{2})+g^{p}(u_{1},u_{2})\right)\\ +\lambda_{N}\int_{-1}^{+\infty}\left(F(x_{1}+b_{1}(x_{1},x_{2})z,x_{2}+b_{2}(x_{1},x_{2})z)-F(x_{1},x_{2})\right)\phi(z)P(z)\mathrm{d}z\end{array}\right\}=0$$
(28)

that is defined for $x_{1},x_{2}\geq 0$. The unique feature of the HJB equation (28) is that there exists a new discount rate $\lambda_{N}\psi^{-1}\int_{0}^{1}D(\phi(z))P(z)\mathrm{d}z$ that depends on the uncertainty, which has not appeared in the literature. This novel term is derived from the discount factor of the Orlicz risk (13).

The so-called Isaacs condition (i.e., the orders of the maximization and minimization are exchangeable) is satisfied in (28). Therefore, we can rewrite (28) as

$$\lambda_Z \inf_{(u_1,u_2)\in C\times C}\left\{F(x_1+c_1(x_1,x_2,u_1),x_2+c_2(x_1,x_2,u_2))-F(x_1,x_2)+g^p(u_1,u_2)\right\}$$

$$+\lambda_N \sup_{\phi(\cdot)>0}\left\{\begin{array}{l}-\psi^{-1}\left(\int_{-1}^{+\infty}D(\phi(z))P(z)\mathrm{d}z\right)F(x_1,x_2)\\+\int_{-1}^{+\infty}\left(F(x_1+b_1(x_1,x_2)z,x_2+b_2(x_1,x_2)z)-F(x_1,x_2)\right)\phi(z)P(z)\mathrm{d}z\end{array}\right\}. \quad (29)$$

$$-\delta F(x_1,x_2)+a_1(x_1,x_2)\frac{\partial F(x_1,x_2)}{\partial x_1}+a_2(x_1,x_2)\frac{\partial F(x_1,x_2)}{\partial x_2}+f^p(x_1,x_2)=0$$

The term inside "sup" is the novel term that arises from the Orlicz risk, which can be further computed as follows. For simplicity, we set

$$\omega(F)(x_1,x_2,z):=\frac{F(x_1+b_1(x_1,x_2)z,x_2+b_2(x_1,x_2)z)-F(x_1,x_2)}{F(x_1,x_2)} \quad (30)$$

if it exists. We assume that $F(x_1,x_2)>0$. Then, we obtain

$$\sup_{\phi(\cdot)>0}\left\{\begin{array}{l}-\psi^{-1}\left(\int_{-1}^{+\infty}D(\phi(z))P(z)\mathrm{d}z\right)F(x_1,x_2)\\+\int_{-1}^{+\infty}\left(F(x_1+b_1(x_1,x_2)z,x_2+b_2(x_1,x_2)z)-F(x_1,x_2)\right)\phi(z)P(z)\mathrm{d}z\end{array}\right\}=-\hat{\delta}(F)(x_1,x_2)F(x_1,x_2), \quad (31)$$

where

$$\hat{\delta}(F)(x_1,x_2):=\frac{\lambda_N}{\psi}\int_{-1}^{+\infty}\left[1-\exp(\psi\omega(F)(x_1,x_2,z))\right]P(z)\mathrm{d}z. \quad (32)$$

The maximizer of (31), which is denoted as $\phi=\hat{\phi}(x_1,x_2,z)$, is expressed as

$$\hat{\phi}(x_1,x_2,z)=\exp(\psi\omega(F)(x_1,x_2,z)). \quad (33)$$

Consequently, equation (29) is rewritten more compactly as

$$\lambda_Z \inf_{(u_1,u_2)\in C\times C}\left\{F(x_1+c_1(x_1,x_2,u_1),x_2+c_2(x_1,x_2,u_2))-F(x_1,x_2)+g^p(u_1,u_2)\right\}$$
$$-\left(\delta+\hat{\delta}(F)(x_1,x_2)\right)F(x_1,x_2)+a_1(x_1,x_2)\frac{\partial F(x_1,x_2)}{\partial x_1}+a_2(x_1,x_2)\frac{\partial F(x_1,x_2)}{\partial x_2}+f^p(x_1,x_2)=0 \quad (34)$$

## 2.6 Another perspective on the HJB equation

The particular form of the HJB equation (28) that possesses the uncertainty-dependent discount rate $\lambda_N\psi^{-1}\int_{-1}^{+\infty}D(\phi(z))P(z)\mathrm{d}z$ suggests that this equation can also be derived as the optimality equation of the following long-run optimal control problem with a control-dependent discount:

$$\text{Find } V(x_1,x_2):=\inf_{\mathbf{u}}\sup_{\phi}\mathbb{E}_{\mathbb{P}(\phi)}\left[\begin{array}{l}\int_0^{+\infty}\exp\left(-\int_0^t(\delta+\tilde{\delta}_s)\mathrm{d}s\right)f^p(X_t^{(1)},X_t^{(2)})\mathrm{d}t\\+\sum_{t\in(0,+\infty),Z}\exp\left(-\int_0^t(\delta+\tilde{\delta}_s)\mathrm{d}s\right)g^p(u_t^{(1)},u_t^{(2)})\end{array}\middle|\mathcal{F}_0\right], \quad (35)$$

with the added novel discount rate that is nonnegative by virtue of $D$:

$$\tilde{\delta}_s=\frac{\lambda_N}{\psi}\int_{-1}^{+\infty}D(\phi_s(z))P(z)\mathrm{d}z\geq 0. \quad (36)$$

Although equation (35) is non-trivial, it is more intuitive than that based on the recursive utility (22) as the former has a more common functional form in control theory. The right-hand side of (35) exists irrespective of $\phi$ because of its nonnegativity and the following upper bound of $V$:

$$V(x_1,x_2) \leq \inf_{\mathbf{u}} \sup_{\phi} \mathbb{E}_{\mathbb{P}(\phi)} \left[ \int_0^{+\infty} \exp\left(-\int_0^t \delta \mathrm{d}s\right) f^p\left(X_t^{(1)}, X_t^{(2)}\right) \mathrm{d}t + \sum_{t \in (0,+\infty), Z} \exp\left(-\int_0^t \delta \mathrm{d}s\right) g^p\left(u_t^{(1)}, u_t^{(2)}\right) \bigg| \mathcal{F}_0 \right] = \frac{1}{\delta} \max_{x_1, x_2 \geq 0} f^p(x_1, x_2) < +\infty \quad (37)$$

This $V$ is clearly nonnegative, and is therefore uniformly bounded.

The optimal controls, which are denoted as $\phi^*, u^{(1),*}, u^{(2),*}$, can be inferred from the HJB equation. Indeed, at each jump time $s$ of process $Z$, we obtain

$$\left(u_s^{(1),*}, u_s^{(2),*}\right) = \arg\min_{(u_1, u_2) \in C \times C} \left\{ F\left(X_s^{(1)} + c_1\left(X_s^{(1)}, X_s^{(2)}, u_1\right), X_s^{(2)} + c_2\left(X_s^{(1)}, X_s^{(2)}, u_2\right)\right) + g^p(u_1, u_2) \right\}, \quad (38)$$

$$\phi_t^* = \hat{\phi}\left(X_t^{(1)}, X_t^{(2)}, z\right), \quad t > 0. \quad (39)$$

Therefore, the analysis of the optimal control is reduced to the resolution of the HJB equation.

## 3. Mathematical analysis and numerical discretization
### 3.1 Exact solution

We study an exact solution to the HJB equation for a single-state case without **Assumption 1**. This simplified case provides useful insights into our stochastic control problem with the unique uncertainty-aversion structure. The SDE can also be viewed as a simplification of jump-driven population dynamics (e.g., [66–68]), whereas our problem concerns the risk and uncertainty.

The HJB equation (34) here is expressed as

$$\lambda_N \sup_{\phi(\cdot)>0} \left\{ -\psi^{-1}\left(\int_{-1}^{+\infty} D(\phi(z))P(z)\mathrm{d}z\right) F(x_1) + \int_{-1}^{+\infty} \left(F(x_1(1+z)) - F(x_1)\right)\phi(z)P(z)\mathrm{d}z \right\}$$
$$-\delta F(x_1) + ax_1 \frac{\partial F(x_1)}{\partial x_1} + (x_1^\alpha)^p = 0$$
, $x_1 > 0$ (40)

with constants $a, \alpha > 0$. The controlled dynamics are

$$\mathrm{d}X_t^{(1)} = X_{t-}^{(1)}\left(a\mathrm{d}t + \mathrm{d}N_t\right), \quad t > 0. \quad (41)$$

Equation (40) is the HJB equation for the linear state dynamics subject to proportional population decay at jump times of process $N$ without human intervention. Proposition 1 demonstrates that equation (40) admits a smooth solution and is identical to (35) under certain conditions. The boundary condition of the solution below is the HJB equation itself.

**Proposition 1** Assume that

$$A := \frac{1}{\delta - \alpha p a - \frac{\lambda_N}{\psi} \int_{-1}^{+\infty} \left(e^{\psi\left\{(1+z)^{\alpha p} - 1\right\}} - 1\right) P(z)\mathrm{d}z} > 0. \quad (42)$$

*Then, equation (40) admits the point-wise solution*

$$F(x_1) = A x_1^{\alpha p}, \quad x_1 > 0, \tag{43}$$

*where $F$ is identical to $V$ in (35). Moreover, the optimal $\phi = \hat{\phi}(x_1, z)$ is independent of $x_1$:*

$$\hat{\phi}(x_1, z) = \exp\left(\psi\left\{(1+z)^{\alpha p} - 1\right\}\right), \quad z > -1. \tag{44}$$

The quantity $F^{1/p}(x_1) = A^{1/p} x_1^{\alpha}$ is the net disutility that is caused by the existence of an exponentially growing population subject to jump perturbations, starting from the initial condition $x_1$. Proposition 2 by direct partial differentiation shows the sensitivity of the solution in Proposition 1.

**Proposition 2** *It follows from (42) that*

$$\frac{\partial A}{\partial \delta} < 0 \quad \text{and} \quad \frac{\partial A}{\partial a} > 0, \tag{45}$$

$$\frac{\partial A}{\partial \lambda_N} \leq 0 \text{ if } P(z) = 0 \text{ for } z > 0 \text{ (resp., } \frac{\partial A}{\partial \lambda_N} \geq 0 \text{ if } P(z) = 0 \text{ for } z < 0\text{),} \tag{46}$$

$$\frac{\partial A}{\partial \psi} \leq 0 \text{ if } P(z) = 0 \text{ for } z > 0 \text{ (resp., } \frac{\partial A}{\partial \psi} \geq 0 \text{ if } P(z) = 0 \text{ for } z < 0\text{).} \tag{47}$$

*Moreover,*

$$\frac{\partial A}{\partial \alpha} > 0 \quad \text{and} \quad \frac{\partial A}{\partial p} > 0 \quad \text{(resp., } \frac{\partial A}{\partial \alpha} < 0 \quad \text{and} \quad \frac{\partial A}{\partial p} < 0 \text{)} \tag{48}$$

*if*

$$0 < a + \lambda_N \int_{-1}^{+\infty} e^{\psi\{(1+z)^{\alpha p}-1\}} (1+z)^{\alpha p} \ln(1+z) P(z) \, dz \quad \text{(resp., } 0 > \text{).} \tag{49}$$

Proposition 2 implies that the disutility increases as the discounting decreases or the growth rate of the population increases according to (45); it increases (resp., decreases) as the intensity of jumps or uncertainty decreases (resp., increases) according to (46)-(47). It follows that the DM should focus sufficiently on the near future ($\delta$ should be sufficiently large) if the speed of the deterministic or stochastic growth is sufficiently high. This assumption is necessary owing to the possibly unbounded dynamics (41), and can be removed under bounded growth, as suggested in (37).

### 3.2 Viscosity property

We demonstrate that the HJB equation (28) is the optimality equation of the differential game (35). For simplicity, we assume the following, which is suggested to be satisfied in our application:

**Assumption 2** *$V$ of (35) is strictly positive for all $x_1, x_2 \geq 0$.*

Assumption 2 implies that any (viscosity) solution of the HJB equation (29) is positive, and the relationship (30) is well defined. The optimality of the HJB equation, where Assumption 2 is not satisfied, was studied in Proposition 1.

**Remark 5** Assumption 2 is satisfied if either $a_1(0,\cdot) > 0$ or $a_1(\cdot,0) > 0$, and $f(x_1,x_2) > 0$ if at least one of $x_1, x_2$ is positive. In this case, the expectation to be optimized in (37) is bounded from below by a positive constant owing to

$$V(x_1,x_2) \geq \inf_{\mathbf{u}} \mathbb{E}_{\mathbb{P}}\left[\begin{array}{l}\int_0^{+\infty} \exp\left(-\int_0^t (\delta+0)\mathrm{d}s\right) f^p\left(X_t^{(1)}, X_t^{(2)}\right)\mathrm{d}t \\ + \sum_{t\in(0,+\infty),Z} \exp\left(-\int_0^t (\delta+0)\mathrm{d}s\right) g^p\left(u_t^{(1)}, u_t^{(2)}\right)\end{array}\bigg| \mathcal{F}_0\right]$$

$$\geq \inf_{\mathbf{u}} \mathbb{E}_{\mathbb{P}}\left[\begin{array}{l}\int_0^{+\infty} \exp(-\delta t) f^p\left(X_t^{(1)}, X_t^{(2)}\right)\mathrm{d}t \\ + \sum_{t\in(0,+\infty),Z} \exp(-\delta t) g^p\left(u_t^{(1)}, u_t^{(2)}\right)\end{array}\bigg| \mathcal{F}_0\right] \quad . \quad (50)$$

$$\geq \inf_{\mathbf{u}} \mathbb{E}_{\mathbb{P}}\left[\int_0^{+\infty} \exp(-\delta t) f^p\left(X_t^{(1)}, X_t^{(2)}\right)\mathrm{d}t \bigg| \mathcal{F}_0\right]$$

The last expectation of (50) never equals zero because it assumes that human intervention is allowed at each Poisson jump times, and there exists a positive probability that the state $\left(X_s^{(1)}, X_s^{(2)}\right)_{s\geq 0}$ belongs to a measurable open set in $(0,+\infty)^2$ for a nonzero time interval (Assumption 1). Thus, the state is immediately pushed to the inside of the domain if it is initially at the boundary ($x_1 x_2 = 0$) owing to $a_1(0,\cdot) > 0$ or $a_1(\cdot,0) > 0$, and such an event should occur with a positive probability. This observation, along with the positivity assumption of $f$, excludes $\inf_{\mathbf{u}} \mathbb{E}_{\mathbb{P}}\left[\int_0^{+\infty} \exp(-\delta t) f^p\left(X_s^{(1)}, X_s^{(2)}\right)\mathrm{d}s \bigg| \mathcal{F}_0\right] = 0$.

The HJB equation satisfies a comparison principle in the sense of constrained viscosity solutions [69–70], such that the supersolution property is satisfied inside and at the boundary of the domain, whereas the subsolution property is required only inside the domain. This asymmetry arises because of the sign condition of the drift and jumps, as stated in Assumption 1. We define the viscosity solutions to our HJB equation based on Definition 4.2 of [69], which are natural solutions to HJB equations that are not necessarily continuously differentiable. See, also [71] for viscosity solutions in general.

We set $\Omega := (0,K)^2$ with a constant $K > 0$ such that $K > \max\{\theta, \bar{b}, \bar{c}\} > 0$ (recall Assumption 1). The closure $\bar{\Omega}$ of $\Omega$ is $\bar{\Omega} = [0,K]^2$. Then, the drift of (1) is always inward along the boundary $\partial\Omega = (\{x_1 = 0, K\} \cap [0,K]) \cup ([0,K] \cap \{x_2 = 0, K\})$. The space of the bounded and uniformly continuous functions in a domain $\omega$ is denoted as $BUC(\omega)$. Similarly, the space of the continuously differentiable functions in $\omega$ is denoted as $C^1(\omega)$. See also [70].

***Definition 1***

*(a) Viscosity supersolution*

Let $\overline{F} \in BUC(\overline{\Omega})$ that is positive on $\overline{\Omega}$. We refer to $\overline{F}$ as a viscosity supersolution of (34) if, for any $(\hat{x}_1, \hat{x}_2) \in \overline{\Omega}$ and $\varphi \in C^1(\overline{\Omega})$ such that $\overline{F}(\hat{x}_1, \hat{x}_2) - \varphi(\hat{x}_1, \hat{x}_2) = \min_{(x_1, x_2) \in \overline{\Omega}} \{\overline{F}(x_1, x_2) - \varphi(x_1, x_2)\}$, it follows that

$$\lambda_Z \inf_{(u_1, u_2) \in C \times C} \{\overline{F}(\hat{x}_1 + c_1(\hat{x}_1, \hat{x}_2, u_1), \hat{x}_2 + c_2(\hat{x}_1, \hat{x}_2, u_2)) - \overline{F}(\hat{x}_1, \hat{x}_2) + g^p(u_1, u_2)\}$$
$$-(\delta + \hat{\delta}(\overline{F})(\hat{x}_1, \hat{x}_2))\overline{F}(\hat{x}_1, \hat{x}_2) + a_1(\hat{x}_1, \hat{x}_2)\frac{\partial \varphi(\hat{x}_1, \hat{x}_2)}{\partial x_1} + a_2(\hat{x}_1, \hat{x}_2)\frac{\partial \varphi(\hat{x}_1, \hat{x}_2)}{\partial x_2} + f^p(\hat{x}_1, \hat{x}_2) \leq 0 \quad (51)$$

*(b) Viscosity subsolution*

Let $\underline{F} \in BUC(\overline{\Omega})$ that is positive on $\overline{\Omega}$. We refer to $\underline{F}$ as a viscosity subsolution of (34) if, for any $(\hat{x}_1, \hat{x}_2) \in \Omega$ and $\varphi \in C^1(\overline{\Omega})$ such that $\underline{F}(\hat{x}_1, \hat{x}_2) - \varphi(\hat{x}_1, \hat{x}_2) = \max_{(x_1, x_2) \in \overline{\Omega}} \{\underline{F}(x_1, x_2) - \varphi(x_1, x_2)\}$, it follows that

$$\lambda_Z \inf_{(u_1, u_2) \in C \times C} \{\underline{F}(\hat{x}_1 + c_1(\hat{x}_1, \hat{x}_2, u_1), \hat{x}_2 + c_2(\hat{x}_1, \hat{x}_2, u_2)) - \underline{F}(\hat{x}_1, \hat{x}_2) + g^p(u_1, u_2)\}$$
$$-(\delta + \hat{\delta}(\underline{F})(\hat{x}_1, \hat{x}_2))\underline{F}(\hat{x}_1, \hat{x}_2) + a_1(\hat{x}_1, \hat{x}_2)\frac{\partial \varphi(\hat{x}_1, \hat{x}_2)}{\partial x_1} + a_2(\hat{x}_1, \hat{x}_2)\frac{\partial \varphi(\hat{x}_1, \hat{x}_2)}{\partial x_2} + f^p(\hat{x}_1, \hat{x}_2) \geq 0 \quad (52)$$

*(c) Viscosity solution*

*A function $F \in BUC(\overline{\Omega})$ is a viscosity solution if it is a viscosity supersolution in terms of (a) and a viscosity subsolution in terms of (b).*

***Remark 6*** The positivity requirement in Definition 1 does not cover Proposition 1, but it is consistent with Remark 5 and the finite difference method presented later. For the exact solution, the boundary condition should alternatively be specified at $x_1 = 0$. We later show that the numerical solutions obtained by our finite- difference method reasonably approximate this boundary condition despite it discretizing the HJB equation over the computational domain.

Hereafter, a viscosity supersolution (resp., viscosity subsolution) is simply referred to as a supersolution (resp., subsolution). Proposition 3 shows the uniqueness of the HJB equation (34).

***Proposition 3*** *There exists at most one viscosity solution in the sense of Definition 1.*

***Remark 7*** The existence of solutions to the HJB equation can be proven, as in Proposition 1, if the solution is smooth. The uniqueness result in Proposition 3 covers both smooth and viscous solutions.

## 3.3 Numerical discretization
### 3.3.1 Finite difference method

As our HJB equation generally cannot be solved analytically, it is discretized by a finite difference method. We prepare a rectangular computational domain $\Omega_{\text{num}} = [0, \bar{X}_1] \times [0, \bar{X}_2]$ with suitably selected constants $\bar{X}_1, \bar{X}_2 > 0$. The grid points in the computational domain $\Omega_{\text{num}}$ are numbered using $N_1, N_2 \in \mathbb{N}$ and the grid sizes $(h_1, h_2) := (\bar{X}_1 / N_1, \bar{X}_2 / N_2)$. In view of Assumption 1, $\bar{X}_1, \bar{X}_2 \geq \max\{\theta, \bar{b}, \bar{c}\}$ should be used, with which no artificial boundary conditions along $x_1 = \bar{X}_1$ and $x_2 = \bar{X}_2$ are required.

We set the vertices $P_{i,j} : (ih_1, jh_2)$ ( $i = 0, 1, 2, ..., N_1$ , $j = 0, 1, 2, ..., N_2$ ) with $(h_1, h_2) := (\bar{X}_1 / N_1, \bar{X}_2 / N_2)$. The solution $F$ discretized at $P_{i,j}$ is denoted as $F_{i,j}$, from which we recover the value function at $P_{i,j}$ as $\Psi_{i,j} = (F_{i,j})^{1/p}$ if $F_{i,j} \geq 0$. The discretization of the HJB equation is carried out using the form of (31), which enables convenient computation. The partial differentials are discretized using the monotone and consistent upwind method (e.g., [72–74]). We assume that the set $C \times C$ is identical to $\Omega_{\text{num}}$ and that $c_i(x_1, x_2, u_1) = -\min\{x_i, u_i\}$ ( $i = 1, 2$ ).

Each term in the HJB equation (29) is discretized as follows (the arrow "$\to$" represents the discretization scheme). Several "$F_{i,j}$" are colored green and are used to update the numerical solution. Several terms are colored magenta for the same reason:

$$-\delta F(x_1, x_2) + f^p(x_1, x_2) \to -\delta F_{i,j} + f^p(ih_1, jh_2), \tag{53}$$

$$a_1(x_1, x_2) \frac{\partial F(x_1, x_2)}{\partial x_1} = \frac{a_1(ih_1, jh_2)}{h_1} \times \begin{cases} F_{i+1,j} - F_{i,j} & (a_1(ih_1, jh_2) \geq 0) \\ F_{i,j} - F_{i-1,j} & (a_1(ih_1, jh_2) < 0) \end{cases}, \tag{54}$$

$$a_2(x_1, x_2) \frac{\partial F(x_1, x_2)}{\partial x_2} = \frac{a_2(ih_1, jh_2)}{h_2} \times \begin{cases} F_{i,j+1} - F_{i,j} & (a_2(ih_1, jh_2) \geq 0) \\ F_{i,j} - F_{i,j-1} & (a_2(ih_1, jh_2) < 0) \end{cases}, \tag{55}$$

$$\lambda_Z \inf_{(u_1, u_2) \in C \times C} \{F(x_1 + c_1(x_1, x_2, u_1), x_2 + c_2(x_1, x_2, u_2)) - F(x_1, x_2) + g^p(u_1, u_2)\}$$

$$= \lambda_Z \inf_{(u_1, u_2) \in C \times C} \{F(x_1 - \min\{x_1, u_1\}, x_2 - \min\{x_2, u_2\}) - F(x_1, x_2) + g^p(u_1, u_2)\}, \tag{56}$$

$$= \lambda_Z \min_{(k,l) \in \{0,1,2,...,N_1\} \times \{0,1,2,...,N_2\}} \{F_{i-\min\{i,k\}, j-\min\{j,l\}} + g^p(kh_1, lh_2)\} - \lambda_Z F_{i,j}$$

and

$$\lambda_N \sup_{\phi(\cdot)>0} \left\{ \begin{array}{l} -\psi^{-1}\left(\int_{-1}^{+\infty} D(\phi(z))p(z)\mathrm{d}z\right)F(x_1,x_2) \\ +\int_{-1}^{+\infty}\left(F(x_1+b_1(x_1,x_2)z, x_2+b_2(x_1,x_2)z)-F(x_1,x_2)\right)\phi(z)p(z)\mathrm{d}z \end{array} \right\}$$

$$= \lambda_N \left\{ \begin{array}{l} -\psi^{-1}\left(\int_{-1}^{+\infty} D(\hat{\phi}(z))p(z)\mathrm{d}z\right)F(x_1,x_2) \\ +\int_{-1}^{+\infty}\left(F(x_1+b_1(x_1,x_2)z, x_2+b_2(x_1,x_2)z)-F(x_1,x_2)\right)\hat{\phi}(z)p(z)\mathrm{d}z \end{array} \right\}\bigg|_{\phi=\hat{\phi}} \quad (57)$$

$$\to \lambda_N \sum_{m=1}^{M} \left\{ -\psi^{-1} D(\hat{\phi}(ih_1, jh_2, z_m)) p_m F_{i,j} + \hat{\phi}(ih_1, jh_2, z_m) p_m (\hat{F}_{i,j,m} - F_{i,j}) \right\}$$

$$= \lambda_N \sum_{m=1}^{M} \hat{\phi}(ih_1, jh_2, z_m) p_m \hat{F}_{i,j,m} - \lambda_N \sum_{m=1}^{M} \left\{ \psi^{-1} D(\hat{\phi}(ih_1, jh_2, z_m)) + \hat{\phi}(ih_1, jh_2, z_m) \right\} p_m F_{i,j}$$

Our method is consistent with Assumption 1 and naturally discretizes the HJB equation considering the direction of characteristics. Given $M \in \mathbb{N}$, we set the discretization $\{z_m\}_{m=1,2,3,\ldots,M}$ of $\mathrm{Supp}P = [\underline{z}, \overline{z}]$ to

$$z_m = \underline{z} + (\overline{z} - \underline{z}) \frac{m - \frac{1}{2}}{M}, \quad m = 1, 2, 3, \ldots, M, \quad (58)$$

and $\hat{F}_{i,j,m}$ is a bilinear interpolation of $F(x_1+b_1(x_1,x_2)z_m, x_2+b_2(x_1,x_2)z_m)$ with vertices $P_{i_m, j_m}$, $P_{i_m+1, j_m}$, $P_{i_m, j_m+1}$, and $P_{i_m+1, j_m+1}$ such that the point $(ih_1+b_1(ih_1,jh_2)z_m, jh_2+b_2(ih_1,jh_2)z_m)$ falls on the inside or boundary of a rectangle that is created by the four vertices.

The optimal controls associated with the (discretized) HJB equation are obtained as follows. The uncertainty is obtained as $\phi(ih_1, jh_2, z_m)$ using (33): for $i = 0,1,2,\ldots,N_1$, $j = 0,1,2,\ldots,N_2$, and $m = 1, 2, 3, \ldots, M$,

$$\hat{\phi}(x_1, x_2, z) \to \hat{\phi}(ih_1, jh_2, z_m) = \exp\left(\psi \frac{\hat{F}_{i,j,m} - F_{i,j}}{\max\{F_{i,j}, \varepsilon\}}\right) \quad (59)$$

with truncation $\max\{\cdot, \varepsilon\}$ ($\varepsilon = 10^{-12}$) to avoid possible division by zero. The discretized control $(u_1, u_2)$ under optimality, which is denoted as $(\hat{u}_1, \hat{u}_2)$, is obtained as follows:

$$\left(\frac{\hat{u}_1}{h_1}, \frac{\hat{u}_2}{h_2}\right) = \arg\min_{(k,l) \in \{0,1,2,\ldots,N_1\} \times \{0,1,2,\ldots,N_2\}} \left\{ F_{i-\min\{i,k\}, j-\min\{j,l\}} - F_{i,j} + g^p(kh_1, lh_2) \right\}. \quad (60)$$

### 3.3.2 Fast sweeping method

The discretized HJB equation is iteratively solved using a fast sweeping method [58,63], which was selected because it is free from the inversion of large matrices, which is usually unavoidable in implicit methods for HJB equations with nonlocal terms. Furthermore, this method is easy to implement because it is basically a fixed-point iteration method. The coefficient multiplied by the green "$F_{i,j}$," which is denoted as $C_{i,j}$, is

$$C_{i,j} = \delta + \lambda_N \sum_{m=1}^{M} \left\{ \psi^{-1} D\left(\hat{\phi}(ih_1, jh_2, z_m)\right) + \hat{\phi}(ih_1, jh_2, z_m) \right\} p_m + \lambda_Z$$
$$+ \frac{a_1(ih_1, jh_2)}{h_1} \times \begin{cases} -1 & (a_1(ih_1, jh_2) \leq 0) \\ 1 & (a_1(ih_1, jh_2) > 0) \end{cases} + \frac{a_2(ih_1, jh_2)}{h_2} \times \begin{cases} -1 & (a_2(ih_1, jh_2) \leq 0) \\ 1 & (a_2(ih_1, jh_2) > 0) \end{cases}. \quad (61)$$

We have $C_{i,j} \geq \delta > 0$. We express the discretized HJB equation at $P_{i,j}$ as

$$F_{i,j} = \frac{H_{i,j}}{C_{i,j}}. \quad (62)$$

Both $H_{i,j}$ and $C_{i,j}$ depend nonlinearly on $F_{i',j'}$ ($i' = 0, 1, 2, ..., N_1$, $j' = 0, 1, 2, ..., N_2$). We use Algorithm A to solve iteratively (62), where $\gamma \in (0,1)$ is a relaxation parameter, the superscript $(n)$ is the iteration count, and $\bar{e} > 0$ is the threshold of iteration errors. Each iteration of Algorithm A (Step 2) is always well defined owing to $C_{i,j} > 0$.

---

*Algorithm A: Fast sweeping method*

**Step 1.** Set an initial guess $\{F_{i,j}^{(0)}\}_{\substack{i=0,1,2,...,N_1 \\ j=0,1,2,...,N_2}}$ and $n = 0$.

**Step 2.** Update $F_{i,j}^{(n+1)} = \gamma \frac{H_{i,j}^{(n)}}{C_{i,j}^{(n)}} + (1-\gamma) F_{i,j}^{(n)}$. The quantities $H_{i,j}^{(n)}$ and $C_{i,j}^{(n)}$ are updated using the Gauss–Seidel method [75]. The sweep direction is from $j = 0$ to $N_2$ (outer loop) and $i = 0$ to $N_1$ (inner loop).

**Step 3.** Evaluate the iteration error $e^{(n)} := \max_{\substack{i=0,1,2,...,N_1 \\ j=0,1,2,...,N_2}} \left| F_{i,j}^{(n+1)} - F_{i,j}^{(n)} \right|$. If $e^{(n)} < \bar{e}$, then terminate the algorithm and output $\{F_{i,j}^{(n+1)}\}_{\substack{i=0,1,2,...,N_1 \\ j=0,1,2,...,N_2}}$ as the numerical solution. Otherwise, set $n \to n+1$ and go to **Step 2**.

---

### 3.3.3 Discussion on discretization

The stability and monotonicity follow from the alternative representation of the discretized HJB equation:

$$G_{i,j}(F) = f^p(ih_1, jh_2), \quad i = 0, 1, 2, ..., N_1, \quad j = 0, 1, 2, ..., N_2, \quad (63)$$

where

$$G_{i,j}(F) := \bar{C}_{i,j}\left(\delta + \lambda_Z + \frac{\lambda_N}{\psi}\sum_{m=1}^{M}\left\{1-\exp\left(-\psi\frac{F_{i,j}-\hat{F}_{i,j,m}}{\max\{F_{i,j},\varepsilon\}}\right)\right\}p_m\right)F_{i,j}$$

$$+\begin{pmatrix}\dfrac{a_1(ih_1,jh_2)}{h_1}\times\begin{cases}F_{i,j}-F_{i+1,j} & (a_1(ih_1,jh_2)\geq 0)\\ F_{i-1,j}-F_{i,j} & (a_1(ih_1,jh_2)<0)\end{cases}\\[2ex] +\dfrac{a_2(ih_1,jh_2)}{h_2}\times\begin{cases}F_{i,j}-F_{i,j+1} & (a_2(ih_1,jh_2)\geq 0)\\ F_{i,j-1}-F_{i,j} & (a_2(ih_1,jh_2)<0)\end{cases}\end{pmatrix}$$

$$-\lambda_Z \min_{(k,l)\in\{0,1,2,\ldots,N_1\}\times\{0,1,2,\ldots,N_2\}}\left\{F_{i-\min\{i,k\},j-\min\{j,l\}}-F_{i,j}+g^p(kh_1,lh_2)\right\}$$

(64)

(some terms are colored green and magenta for explanation purposes). It follows that $G_{i,j}(F)$ is increasing with respect to the green "$F_{i,j}$" if it is positive, but it is nontrivial a priori (i.e., the function $y(1-e^{1/y})$ is non-decreasing for all $y>0$). It is also increasing with respect to the magenta "$F_{i,j}-F_{\cdot,\cdot}$"s if all $F_{i,j}$ are positive. We assume that this nonnegativity property is satisfied. Then, combined with Remark 3.10, Lemma 3.11 of [76] shows that the proposed finite difference method is nonnegative and monotone in the sense of Definition 3.9 if all $F_{i,j}$ are positive. Its stability, namely the uniform boundedness with respect to the discretization parameter, follows (e.g., repeat the argument in the proof of Lemma 5.1 of [77]):

$$0 \leq \min_{(i,j)\in\{0,1,2,\ldots,N_1\}\times\{0,1,2,\ldots,N_2\}} f^p(ih_1,jh_2) \leq \delta F_{i,j} \leq \max_{(i,j)\in\{0,1,2,\ldots,N_1\}\times\{0,1,2,\ldots,N_2\}} f^p(ih_1,jh_2). \quad (65)$$

The left-most inequality arises from the nonnegativity of $f$.

Proposition 4 below is the main contribution of this study from a theoretical numerical analysis perspective. The proof is based on a modified equation with truncated coefficients, which is stable and monotone, and further admits a unique numerical solution (e.g., [78]). We demonstrate that this solution is desired and satisfies (63).

***Proposition 4*** *The proposed finite difference method that employs the discretization (63) satisfies a unique solution that complies with (65).*

***Remark 8*** The consistency of the finite difference method is also important theoretically. This is complicated in the proposed HJB equation because of the boundary treatment. We expect that it holds if the solution to the HJB equation is smooth.

## 4. Application
### 4.1 Comparison with exact solution
The numerical solutions that are generated by the finite difference method are examined against the exact solution. The HJB equation is the more realistic logistic version of (40), where the domain is the bounded interval $(0,1)$ and the drift is nonlinear, and the jump perturbation decreases the population ($P(z)=1$ for $z \in (-1, 0)$ and $P(z) = 0$ otherwise) with $\lambda_N = 0.05$:

$$\lambda_N \sup_{\phi(\cdot) > 0} \left\{ -\psi^{-1} \left( \int_{-1}^{+\infty} D(\phi(z)) P(z) \mathrm{d}z \right) F(x_1) + \int_{-1}^{+\infty} \left( F(x_1(1+z)) - F(x_1) \right) \phi(z) P(z) \mathrm{d}z \right\} \\ -\delta F(x_1) + ax_1 (1 - x_1) \frac{\mathrm{d} F(x_1)}{\mathrm{d} x_1} + (x_1^\alpha)^p = 0$$
(66)

The solution to (66) does not exactly coincide with that in Proposition 1 because the drift coefficients differ. It is expected that they asymptotically agree well for a sufficiently small $x_1 > 0$ because $ax_1(1 - x_1)$ can be approximated as $ax_1$ for such $x_1 > 0$. Therefore, the asymptotic behavior of the numerical solutions can be used to evaluate the performance of the finite difference method.

We use the parameter values $\delta = 0.1$, $a = 0.02$, $p = 2$, and $\psi = 0.01$ for the HJB equation, and $\gamma = 0.5$ and $\bar{e} = 10^{-10}$ for Algorithm 1. The computational resolutions are $M = N_1 = 100$ and $M = N_1 = 500$. We examine three qualitatively different cases of the solution profile considering Proposition 1: the convex case ($p\alpha = 2$), linear case ($p\alpha = 1$), and concave case ($p\alpha = 0.5$). **Figure 1** compares the numerical solutions and exact solution of Proposition 1, whereas **Figure 2** magnifies the plots in **Figure 1** near the origin $x_1 = 0$. **Figure 1** indicates that the numerical and exact solutions agree well near the origin, whereas they deviate for larger $x_1$ values, as expected, owing to the emergence of the nonlinear drift in the numerical solutions. **Figure 2** confirms the agreement between the numerical and exact solutions for small values, and further demonstrates that the numerical solutions are almost identical for the two computational resolutions. Therefore, the computational resolution of this order (i.e., $M = N_1 = O(10^2)$) can be satisfactory. The total iteration counts $n$ for obtaining the converged numerical solutions in the concave case are 179 for $M = N_1 = 100$ and 633 for $M = N_1 = 500$. Similar results are obtained for the linear and concave cases (convex case: 175 and 613, linear case: 178 and 626).

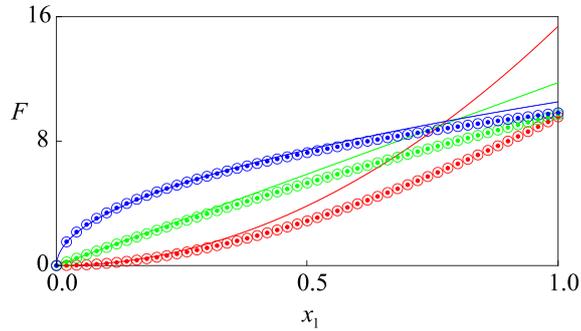

**Figure 1.** Comparison of exact (curves) and numerical (circles) solutions for convex (red), linear (green), and concave (blue) cases. The filled and unfilled circles represent computational resolutions of $M = N_1 = 100$ and $M = N_1 = 500$, respectively.

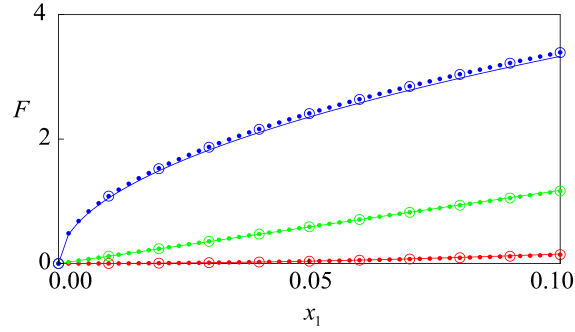

**Figure 2.** Magnification of results in **Figure 1** near origin. Note that the exact and numerical solutions are only asymptotically the same.

## 4.2 Aquatic vegetation management
### 4.2.1 Study area

The study area for the model application was Lake Shinji, which is located at 35.4426N and 132.9514E in the estuary of the Hii River System along the Sea of Japan (**Figure 3**). Lake Shinji is a brackish lake with an area of 79.1 (km$^2$), shoreline length of 47 (km), and mean depth of 4.5 (m), and its basin has the largest drainage area of 1,288.4 (km$^2$) among the brackish lakes in southwestern Japan. This lake has been registered as a Wetland of International Importance under the Ramsar Convention [79], and it has been studied to analyze and monitor environmental and ecosystem states (e.g., [80–81]).

The submerged aquatic vegetation *Potamogeton anguillanus* is a representative macrophyte of Lake Shinji. *P. anguillanus* was distributed in remote patches in Lake Shinji until the 1990s, and then began proliferating around 2010 [82]. It has recently been overgrowing throughout the lake annually. Owing to its hybrid origin, this macrophyte generally reproduces asexually through various methods, including vegetative propagation, thalli fragmentation, and rarely, turion sprout and seed apomixis. Vegetative propagation and thalli fragmentation involve the growth of new macrophytes from rhizomes and fragmented thalli, respectively. This macrophyte is often found in nutrient-rich environments and its overgrowth reduces light penetration and oxygen levels, which leads to the mortality of other phytes and algae as well as interference with inland fisheries and recreational activities in Lake Shinji. Since 2020, we have conducted periodic ecological surveillance of *P. anguillanus* using visual and aerial observations at seven distant stations along the shores (**Figure 4**). Moreover, both the national and local governments have been seeking methods to suppress the overgrowth in a cost-efficient manner.

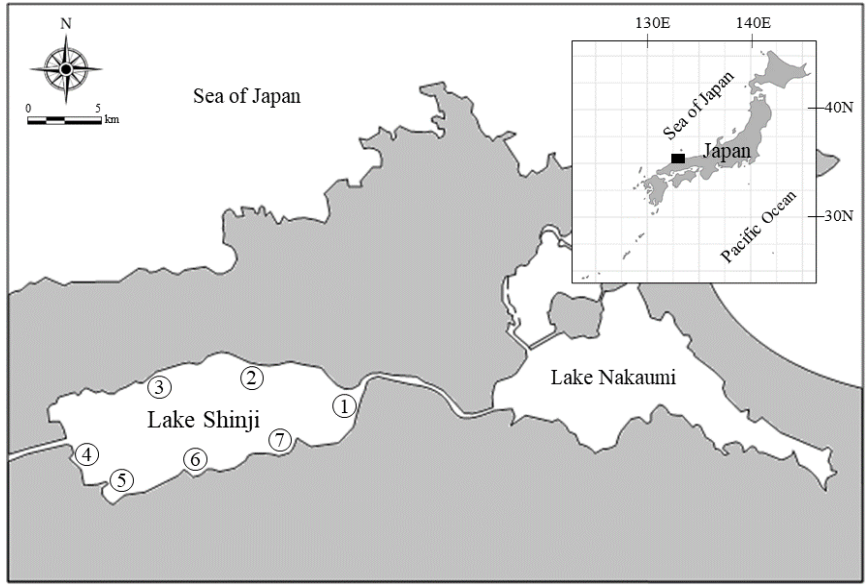

**Figure 3.** Map of seven ecological surveillance stations (①　to ⑦) of *P. anguillanus* in Lake Shinji, western Japan.

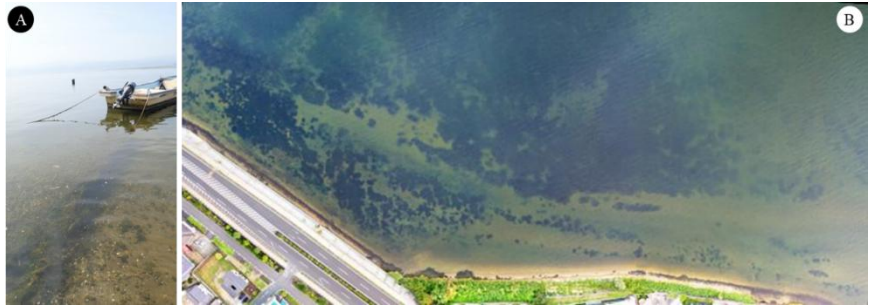

**Figure 4.** (A) Visual observation of *P. anguillanus* overgrowth interfering with boat navigation at Station 2 and (B) aerial panoramic observation of huge deep-greenish colony of *P. anguillanus* at Station 7 in Lake Shinji.

### 4.2.2 Parameter estimation

We mathematically describe the population dynamics of *P. anguillanus* at each station considering the growing (superscript "1") and drifting (superscript "2") populations. The system dynamics are based on **Remark 1**, with an additional source term for the growing population:

$$a_1(x_1, x_2) = r_1 x_1 \left(1 - \frac{x_1}{Q_1}\right) - d_1 x_1 + \alpha_1 \max\{Q_1 - x_1, 0\}, \quad a_2(x_1, x_2) = \chi_\mu(x_2)\left\{r_2 x_2 \left(1 - \frac{x_2}{Q_2}\right) + d_1 x_1\right\}, \qquad (67)$$

with $b_i$ in (6) and $c_i$ in (7). We assume the equalities $Q_1 = Q_2$ and $r_1 = r_2$ owing to the lack of data on the population dynamics of the drifting vegetation. Therefore, the parameters to be fitted are $r_1$, $d_1$, $\alpha_1$, and $Q_1$. The subscript "1" is omitted below to simplify the description.

We estimated the model parameters based on the field survey results obtained in 2022. The growing macrophyte population at each station was estimated by calculating the covering ratio using the belt-transect method, such that an $120 \times 4 = 480$ (m²) rectangular region along the shoreline was divided into 480 squares with an area of $1 \times 1 = 1$ (m²). The survey was conducted from April to September in 2022, and we evaluated the population of growing macrophytes based on the number of squares in which the macrophyte populations were visible (see Appendix B). At each station, the parameters were fitted against the forward Euler discretization of the deterministic model:

$$X(t_{i+1}) = \max\left\{0, X(t_i) + (t_{i+1} - t_i)\left(rX(t_i)\left(1 - \frac{X(t_i)}{Q}\right) - dX(t_i) + \alpha \max\{Q - X(t_i), 0\}\right)\right\}, \quad (68)$$

with $t_i$ for the $i$ th observation day ($i = 1, 2, 3, ..., I$), $I$ is the total number of observations, and $X(t_i)$ is the population of the growing vegetation on that day. The observation value is denoted with the superscript "obs" as $X^{obs}$. Assuming that $X(t_1) = X^{obs}(t_1)$, we used the least-squares method to minimize $\sum_{i=1}^{I}(X(t_i) - X^{obs}(t_i))^2$.

**Table 1** summarizes the obtained parameter values for all stations (see also Appendix B). The identification results demonstrate variability in the parameter values, and hence, the growth dynamics among the surveyed stations. The growth rate $r$ and decay rate $d$ vary by three orders, whereas the capacity $Q$ is almost 1 except for Station 3. The external source $\alpha$ varies by two orders among the stations. It has been suggested that Stations 3 and 7 are particularly important sites for fisheries in Lake Shinji, but they are often subject to the thick growth of the macrophytes that disturb the near-shore ship sailing for fisheries activities. Indeed, **Table 1** suggests a relatively high growth rate $r$ and transport rate $d$ at Stations 3 and 7. Finally, the estimation procedure of the model does not consider jump perturbation. Such jumps occur in typhoon seasons during which the vegetation patches are subject to extreme hydrodynamic disturbance.

Table 1. Parameter values for Stations 1–7.

|  | $r$ (1/day) | $d$ (1/day) | $\alpha$ (1/day) | $Q/480$ (-) |
|---|---|---|---|---|
| Station 1 | 4.88.E-03 | 4.11.E-04 | 4.19.E-04 | 9.88.E-01 |
| Station 2 | 1.78.E-02 | 2.72.E-02 | 1.49.E-03 | 9.83.E-01 |
| Station 3 | 1.14.E+00 | 5.05.E-01 | 1.31.E-02 | 6.71.E-01 |
| Station 4 | 1.91.E-02 | 3.22.E-04 | 5.01.E-04 | 9.95.E-01 |
| Station 5 | 7.91.E-02 | 6.15.E-02 | 9.30.E-04 | 9.88.E-01 |
| Station 6 | 1.85.E-01 | 1.29.E-01 | 1.29.E-03 | 9.17.E-01 |
| Station 7 | 4.01.E-01 | 2.31.E-01 | 1.73.E-04 | 9.98.E-01 |

### 4.2.3 Numerical computation

We first analyzed all seven stations to demonstrate that the optimal harvesting policies differed significantly, even under the same objective function. Subsequently, we conducted a sensitivity analysis on Station 7 as a key station from a fisheries perspective.

The domain was the square $(0,1)^2$ with normalization by $K = 480$ (m²), and was uniformly discretized with 301 vertices in each direction (see, Appendix C). We focused on the controls $(u_1, u_2)$ and added discount rate (36). Unless otherwise specified, the parameter values were as follows, which were selected to obtain nontrivial (i.e., nonzero) optimal controls are obtained (e.g., "do nothing" $(u_1, u_2) = (0, 0)$ may become optimal if the intervention cost is extremely large): $\lambda_Z = 1/30$ (day), $\lambda_N = 1/50$ (day), $\delta = 1/30$ (day), $p = 2$, and $\psi = 10$. The parameter values for the population dynamics were set according to **Table 1**. The coefficient $f$ was set to $f(x_1, x_2) = |x_1 - 1/4| + x_2$ assuming that exterminating the macrophytes in the lake is not optimal (first term) and the drifted population should be cleaned up to reduce the fisheries disutility (second term). The probability density function $P$ was set to $2(1-z)$ in $[0,1]$, which is the simplest model in which smaller disturbances more likely occur. The coefficient $g$ was set to $g(u_1, u_2) = c_0 + c_1(u_1 + u_2)$, with the fixed observation cost $c_0 = 1$ incurred at each observation and the harvesting cost $c_1(u_1 + u_2)$ with $c_1 = 5$ proportional to the harvested population. This proportionality is justified because the harvesting effort would be proportional to the area to be harvested.

**Figures 5**, **6**, and **7** depict the computed optimal $u_1$, $u_2$, and added discount rate $\bar{\delta}$, respectively, at each station. The white area denotes $u_1 = 0$ (do nothing against the growing population). **Figures 5** and **6** illustrate that the optimal harvesting policies for the macrophytes differ significantly among the seven stations. As can be observed from **Figure 5**, for a growing population with a relatively small growth rate $r$ (Stations 1, 2, and 4), it is optimal to harvest the population only if it is sufficiently large, although the free boundary to activate the harvesting at each observation is complex and does not appear to be simply parameterized. Nevertheless, the DM can decide on the harvesting policy as a function of the growing and drifting populations $(x_1, x_2)$ at each observation time. Station 3 differs from the others because "do nothing" is optimal, which is owing to the extremely high growth rate among the stations for which harvesting the growing population becomes costly. For Stations 6 and 7 with smaller but relatively large growth rates, the growing population should be harvested if it is smaller than a certain threshold Station 5 is also unique because harvesting the growing population is always optimal, whereas the amount of harvesting depends on the state. This type of difference between stations cannot be obtained without solving the HJB equation. **Figure 6** indicates that the optimal harvesting policy against the drifting population is complex for all stations, although there are some similarities among stations; for example, for stations 1, 2, and 4, it is optimal to do nothing if both the growing and drifting populations are small, which is reasonable considering the functional shape of $f$. Another reason for this finding is that these stations have relatively small growth rates compared to the others. **Figure 7** demonstrates that the added discount rate $\bar{\delta}$ does not exhibit abrupt changes in the computational domain; hence, it is inferred to be a continuous function at all stations. It has a non-monotone profile and is difficult to parameterize in a straightforward manner. Nevertheless, a common result is that the added discount has an extreme point in the region where "do nothing" against both the growing and drifting populations is optimal. The determination of the optimal harvesting policy without solving the HJB equation would not be effective because of the complex dependence of the optimal controls on the state.

Finally, we conducted a sensitivity analysis for Station 7, as illustrated in **Figures 8** and **9**. We examined (a) the nominal case, which is the same as those in **Figures 6** and **7**, (b) stronger uncertainty aversion with $\psi = 10^{3/2}$, (c) a smaller cost with $c_0 = 1/2$ and $c_1 = 5/2$, (d) a larger discount with $\delta = 1/15$ (1/day), (e) more frequent observation with $\lambda_Z = 1/15$ (1/day), and (f) stronger risk aversion with $p = 2.5$. The free boundary, which is independent of the drifting population, manifests in different locations in each case for the optimal harvesting policy against the growing population. **Figure 8** suggests that the conditions under which the growing population should be harvested become wider if the harvesting cost is reduced or more frequent observations are possible, which implies that technological progress leads to a more active management policy against the growing population. In contrast, stronger risk or uncertainty aversion or a more myopic mindset leads to less active optimal harvesting of the growing population. According to **Figure 9**, less active management against the growing population is associated with a simpler form of harvesting against the drifting population, which suggests that it is ineffective to plan a complex management policy under large risk or uncertainty. This result is considered due to the lower future predictability of the system to be controlled by the DM, which suggests that a sophisticated management policy should be implemented under sufficiently small risk or uncertainty. Nevertheless, such preliminary activities would incur additional management costs, which are beyond the scope of this study but would be an interesting topic from theoretical and practical perspectives.

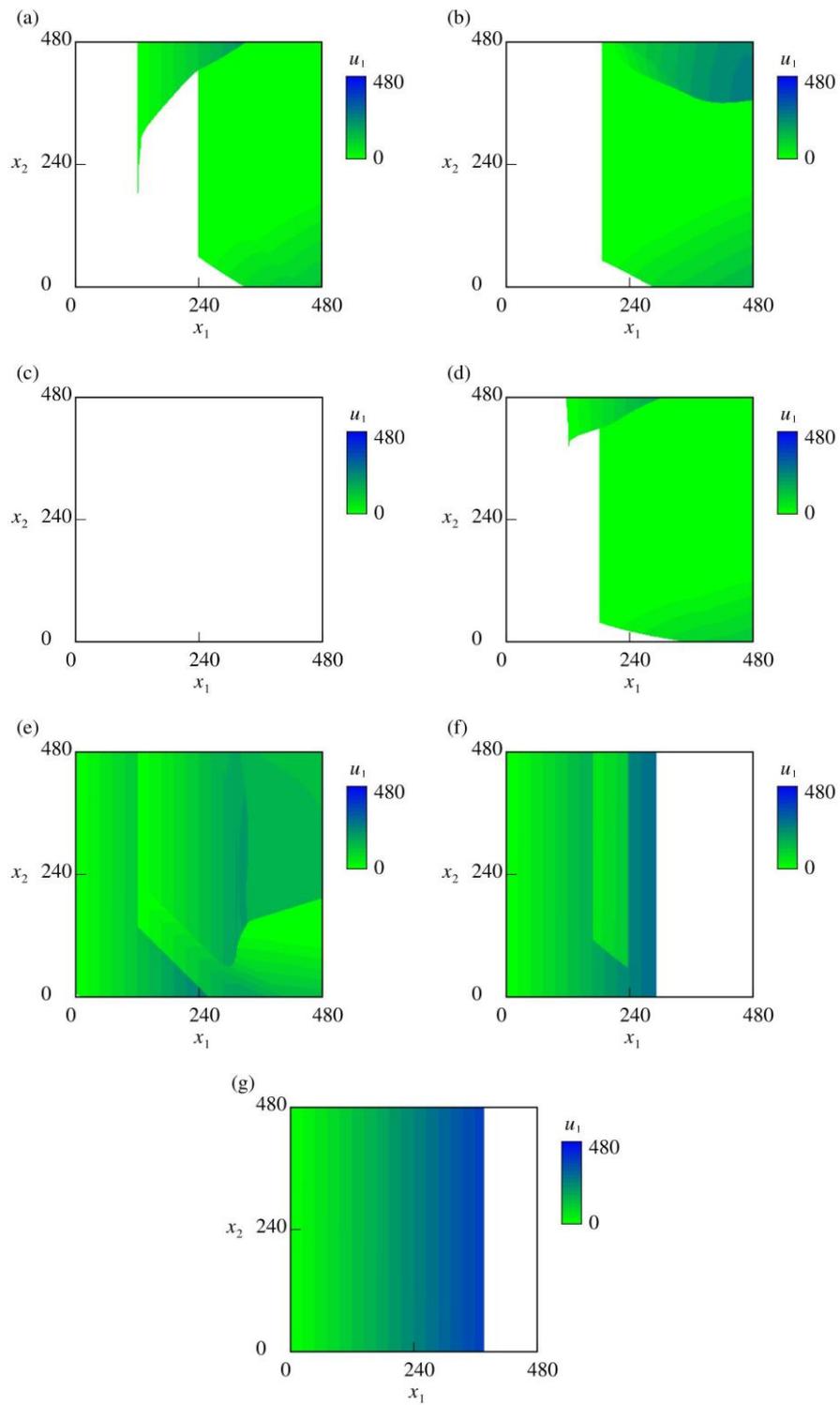

**Figure 5.** The computed optimal $u_1$ at Stations (a) 1, (b) 2, (c) 3, (d) 4, (e) 5, (f) 6, and (g) 7. The white area represents $u_1 = 0$ (do nothing against the growing population).

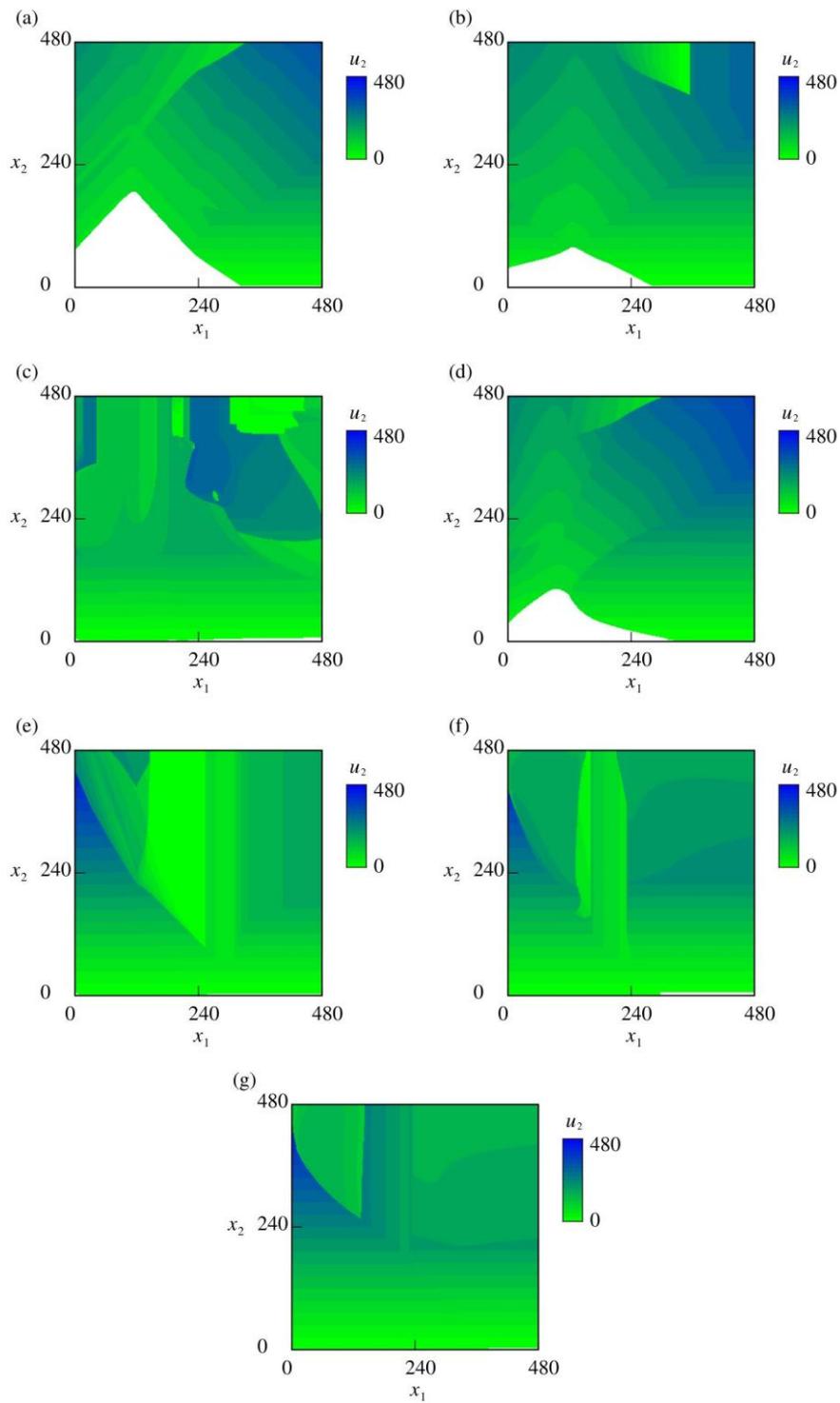

**Figure 6.** Computed optimal $u_2$ at Stations (a) 1, (b) 2, (c) 3, (d) 4, (e) 5, (f) 6, and (g) 7. The white area represents $u_2 = 0$ (do nothing against the drifting population).

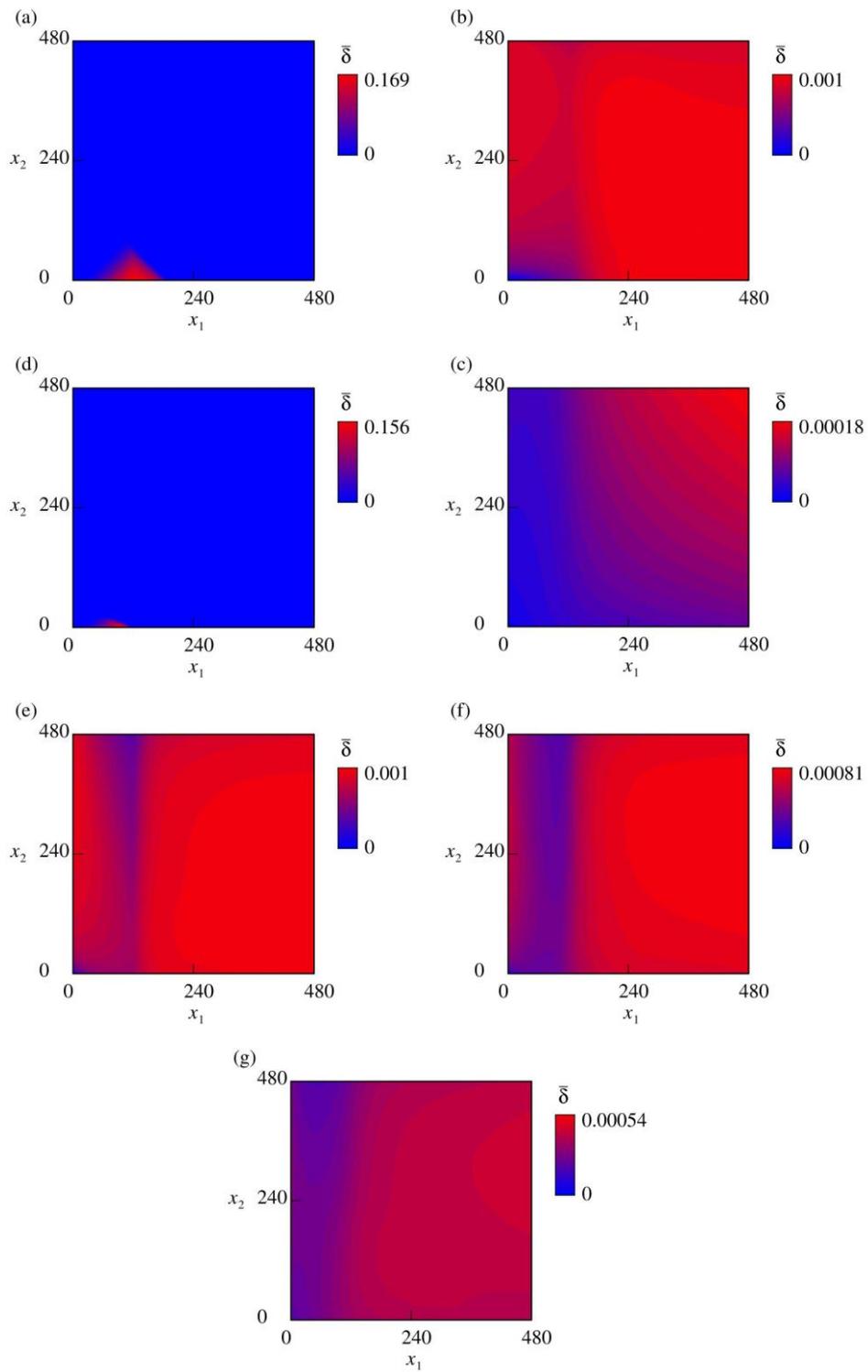

**Figure 7.** Added discount rate $\bar{\delta}$ at Stations (a) 1, (b) 2, (c) 3, (d) 4, (e) 5, (f) 6, and (g) 7.

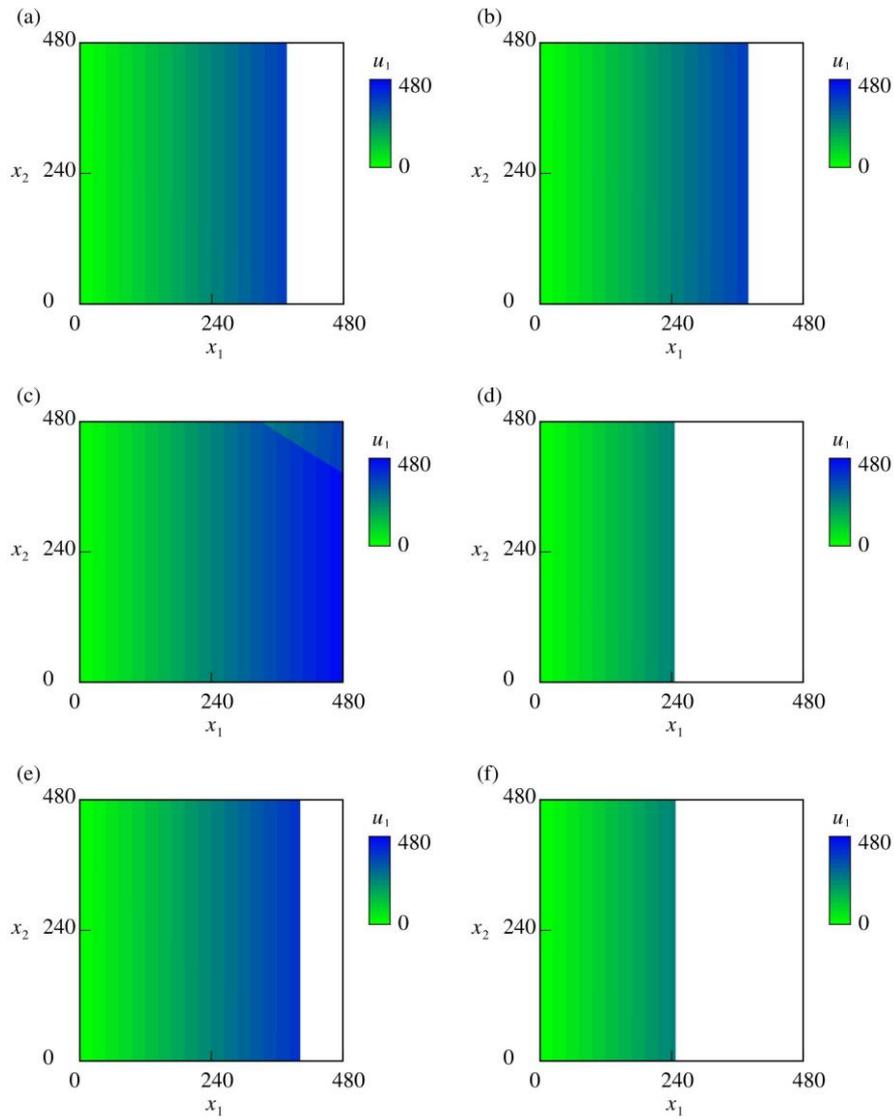

**Figure 8.** Computed optimal $u_1$ at Station 7 under different conditions: (a) nominal, (b) stronger uncertainty aversion, (c) smaller cost, (d) larger discount, (e) more frequent observation, and (f) stronger risk aversion. The white area represents $u_1 = 0$ (do nothing against the growing population).

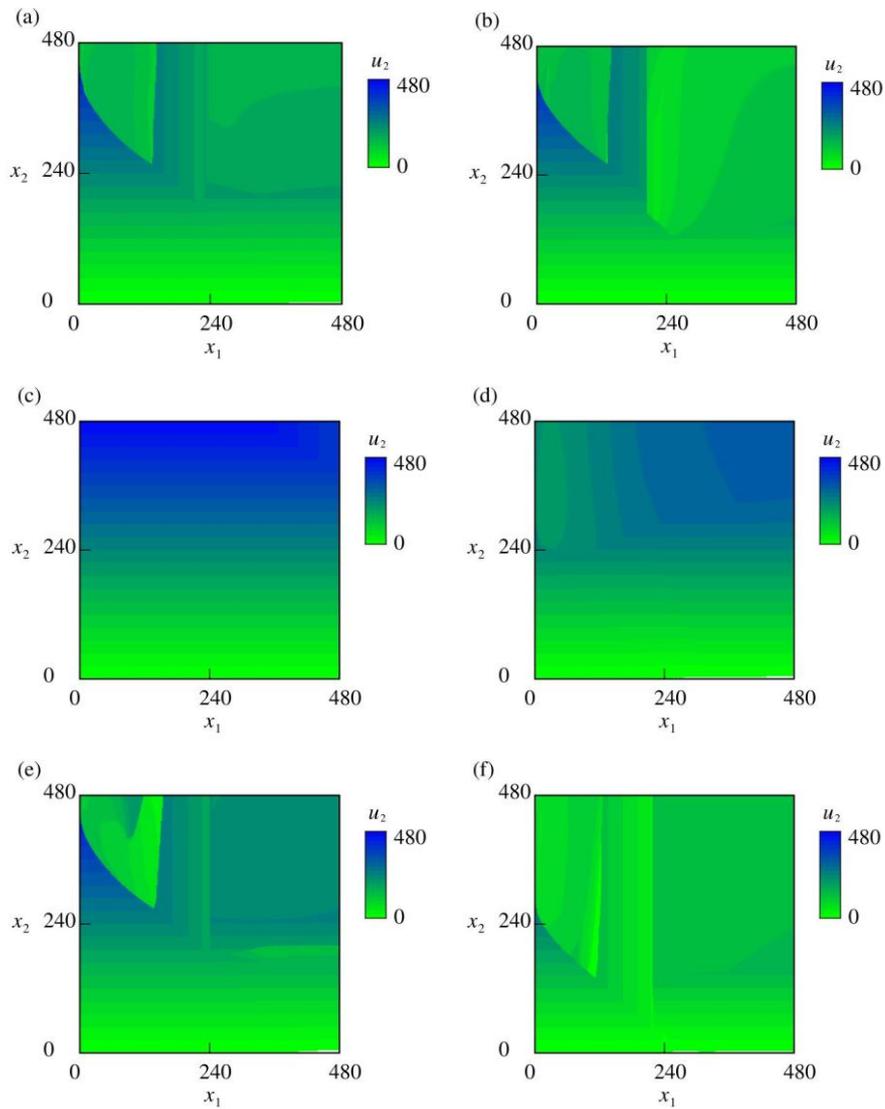

**Figure 9.** Computed optimal $u_2$ at Station 7 under different conditions: (a) nominal, (b) stronger uncertainty aversion, (c) smaller cost, (d) larger discount, (e) more frequent observation, and (f) stronger risk aversion. The white area represents $u_2 = 0$ (do nothing against the growing population).

## 5. Conclusions

The robustified dynamic Orlicz risk was applied to the stochastic control problem in aquatic environmental management and restoration subject to risk and uncertainty. The hybrid SDE that represents the state dynamics was controlled through impulsive interventions. The HJB equation with a unique state-dependent discount rate was derived from the discrete-time formulation of the robustified dynamic Orlicz risk. The uniqueness of this HJB equation was discussed and a finite difference method for its numerical discretization was presented. Furthermore, its convergence was verified using the exactly solvable case. Finally, the management of macrophytes in a lake was investigated using the HJB equation.

The proposed mathematical approach can be applied to other problems in which sustainability is a major concern, such as tourism [83], fruit harvesting [84], fisheries management [85], social discounting under climate changes [86], and inventory [87]. Although our analysis focused on the bounded state dynamics, the unbounded dynamics would also naturally arise in applications. The Young functions should be selected more carefully in such cases because both the risk and uncertainty are sensitive to the far-field behavior of the state dynamics [88]. Other divergences [89–90] can be used instead of the Kullback–Leibler divergence, but the well-posedness of the problem in such cases is nontrivial. In our problem, the stability was theoretically and computationally guaranteed owing to an exogenous positive discount. The vanishing discount limit is also of interest in applications as it corresponds to the time-average control [91] that can be more intuitively connected to the concept of sustainability. In the future, this case will be investigated based on our framework.

Interesting research topics in other areas include the preventive maintenance of industrial machinery [92], grazing ecosystem dynamics [93], and epidemic spread [94], in which jump-driven dynamics play a role. Theoretically, the HJB equation in the extended case is a highly nontrivial equation, and the numerical discretization of such an equation can be computationally challenging. In future work, we will address these issues by focusing on long-term management problems in river environments.

# Appendices of "Environmental management and restoration under unified risk and uncertainty using robustified dynamic Orlicz risk" by Hidekazu Yoshioka, Motoh Tsujimura, Futoshi Aranishi, and Tomomi Tanaka

**Appendix A: Proofs of Propositions**

*Proof of Proposition 1*

The equation (40) is rewritten as

$$-\frac{\lambda_N}{\psi}\int_{-1}^{+\infty}\left(1-\exp\left(\psi\frac{F(x_1(1+z))-F(x_1)}{F(x_1)}\right)\right)P(z)\mathrm{d}z\times F(x_1)-\delta F(x_1)$$
$$+ax_1\frac{\mathrm{d}F(x_1)}{\mathrm{d}x_1}+(x_1^\alpha)^p=0 \quad , \quad x_1>0. \tag{69}$$

Substituting (43) to (69) yields

$$-\frac{\lambda_N}{\psi}\int_{-1}^{+\infty}\left(1-\exp\left(\psi\frac{A(x_1(1+z))^{\alpha p}-A(x_1)^{\alpha p}}{A(x_1)^{\alpha p}}\right)\right)P(z)\mathrm{d}z\times A(x_1)^{\alpha p}-\delta A(x_1)^{\alpha p}$$
$$+ax_1\times\alpha pA(x_1)^{\alpha p-1}+(x_1^\alpha)^p=0 \quad , \quad x_1>0. \tag{70}$$

This equation is rewritten as

$$-\frac{\lambda_N}{\psi}\int_{-1}^{+\infty}\left(1-\exp\left(\psi\left\{(1+z)^{\alpha p}-1\right\}\right)\right)P(z)\mathrm{d}z\times A-\delta A+\alpha paA+1=0, \tag{71}$$

or equivalently

$$\left(\delta-\alpha pa-\frac{\lambda_N}{\psi}\int_{-1}^{+\infty}\left(e^{\psi\{(1+z)^{\alpha p}-1\}}-1\right)P(z)\mathrm{d}z\right)A=1. \tag{72}$$

Rearranging (72) yields (42). Because the solution itself is smooth and grows polynomially, the fact that it is identical to the optimized expectation in (35) can be checked following the standard verification argument (Chapter 9.1 of Øksendal and Sulem [14]). We check that the right-hand side of (35) exists. To see this, notice that the controlled dynamics (41) read

$$X_t^{(1)}=x_1\exp\left(at+\sum_{s\in(0,t),N}\ln(1+\Delta N_s)\right), \quad t>0. \tag{73}$$

We also have

$$\tilde{\delta}_s\Big|_{\phi=\hat{\phi}}=\frac{\lambda_N}{\psi}\int_{-1}^{+\infty}D(\hat{\phi}(z))P(z)\mathrm{d}z$$
$$=\frac{\lambda_N}{\psi}\int_{-1}^{+\infty}\left(\psi\left\{(1+z)^{\alpha p}-1\right\}e^{\psi\{(1+z)^{\alpha p}-1\}}-e^{\psi\{(1+z)^{\alpha p}-1\}}+1\right)P(z)\mathrm{d}z. \tag{74}$$
$$(:=\overline{\delta})$$

Thereby, we obtain

$$\mathbb{E}_{\mathbb{P}(\hat{\phi})}\left[\begin{array}{l}\int_0^{+\infty}\exp\left(-\int_0^t(\delta+\tilde{\delta}_s)\mathrm{d}s\right)f^p\left(X_t^{(1)},X_t^{(2)}\right)\mathrm{d}t\\+\sum_{t\in(0,+\infty),Z}\exp\left(-\int_0^t(\delta+\tilde{\delta}_s)\mathrm{d}s\right)g^p\left(u_t^{(1)},u_t^{(2)}\right)\end{array}\bigg|\mathcal{F}_0\right]$$

$$=\mathbb{E}_{\mathbb{P}(\hat{\phi})}\left[\int_0^{+\infty}\exp\left(-\int_0^t(\delta+\tilde{\delta}_s)\mathrm{d}s\right)\left(X_t^{(1)}\right)^{\alpha p}\mathrm{d}t\bigg|\mathcal{F}_0\right]$$

$$=\mathbb{E}_{\mathbb{P}(\hat{\phi})}\left[\int_0^{+\infty}\exp\left(-(\delta+\bar{\delta})t\right)\left(X_s^{(1)}\right)^{\alpha p}\mathrm{d}t\bigg|\mathcal{F}_0\right] \qquad (75)$$

$$=(x_1)^{\alpha p}\mathbb{E}_{\mathbb{P}(\hat{\phi})}\left[\int_0^{+\infty}\exp\left(-(\delta+\bar{\delta}-a\alpha p)t\right)\exp\left(\alpha p\sum_{s\in(0,t),N}\ln(1+\Delta N_s)\right)\mathrm{d}t\bigg|\mathcal{F}_0\right]$$

$$=(x_1)^{\alpha p}\mathbb{E}_{\mathbb{P}(\hat{\phi})}\left[\int_0^{+\infty}\exp\left(-(\delta+\bar{\delta}-a\alpha p)t\right)\exp\left(\sum_{s\in(0,t),N}\ln(1+\Delta N_s)^{\alpha p}\right)\mathrm{d}t\bigg|\mathcal{F}_0\right]$$

$$=(x_1)^{\alpha p}\int_0^{+\infty}\exp\left(-(\delta+\bar{\delta}-a\alpha p)t\right)\mathbb{E}_{\mathbb{P}(\hat{\phi})}\left[\exp\left(\sum_{s\in(0,t),N}\ln(1+\Delta N_s)^{\alpha p}\right)\right]\mathrm{d}t$$

We have

$$\mathbb{E}_{\mathbb{P}(\hat{\phi})}\left[\exp\left(\sum_{s\in(0,t),N}\ln(1+\Delta N_s)^{\alpha p}\right)\right]=\exp\left(\left(\lambda_N\int_{-1}^{+\infty}\left\{(1+z)^{\alpha p}-1\right\}e^{\psi\left\{(1+z)^{\alpha p}-1\right\}}P(z)\mathrm{d}z\right)t\right). \qquad (76)$$

Then, we obtain

$$\mathbb{E}_{\mathbb{P}(\hat{\phi})}\left[\begin{array}{l}\int_0^{+\infty}\exp\left(-\int_0^t(\delta+\tilde{\delta}_s)\mathrm{d}s\right)f^p\left(X_t^{(1)},X_t^{(2)}\right)\mathrm{d}t\\+\sum_{t\in(0,+\infty),Z}\exp\left(-\int_0^t(\delta+\tilde{\delta}_s)\mathrm{d}s\right)g^p\left(u_t^{(1)},u_t^{(2)}\right)\end{array}\bigg|\mathcal{F}_0\right], \qquad (77)$$

$$=(x_1)^{\alpha p}\int_0^{+\infty}\exp\left(-\left(\delta+\bar{\delta}-a\alpha p-\lambda_N\int_{-1}^{+\infty}\left\{(1+z)^{\alpha p}-1\right\}e^{\psi\left\{(1+z)^{\alpha p}-1\right\}}P(z)\mathrm{d}z\right)t\right)\mathrm{d}t$$

where

$$\delta+\bar{\delta}-a\alpha p-\lambda_N\int_{-1}^{+\infty}\left\{(1+z)^{\alpha p}-1\right\}e^{\psi\left\{(1+z)^{\alpha p}-1\right\}}P(z)\mathrm{d}z$$

$$=\delta-a\alpha p+\frac{\lambda_N}{\psi}\int_{-1}^{+\infty}\left(\psi\left\{(1+z)^{\alpha p}-1\right\}e^{\psi\left\{(1+z)^{\alpha p}-1\right\}}-e^{\psi\left\{(1+z)^{\alpha p}-1\right\}}+1\right)P(z)\mathrm{d}z$$

$$-\frac{\lambda_N}{\psi}\int_{-1}^{+\infty}\psi\left\{(1+z)^{\alpha p}-1\right\}e^{\psi\left\{(1+z)^{\alpha p}-1\right\}}P(z)\mathrm{d}z \qquad (78)$$

$$=\delta-a\alpha p-\frac{\lambda_N}{\psi}\int_{-1}^{+\infty}\left(e^{\psi\left\{(1+z)^{\alpha p}-1\right\}}-1\right)P(z)\mathrm{d}z$$

$$=\frac{1}{A}$$

$$>0$$

Consequently, we arrive at the desired result due to

$$V(x_1)=(x_1)^{\alpha p}\int_0^{+\infty}\exp\left(-\frac{1}{A}t\right)\mathrm{d}t=A(x_1)^{\alpha p}. \qquad (79)$$

□

*Proof of Proposition 3*

The proof follows that of Theorem 5.1 of Tian and Ba [69] (this literature is simply called TB22 in this proof for simplicity). The difference between our and their HJB equations is that the former has a positive discount rate ($\delta > 0$), while the latter does not. The existence of the positive discount rate can avoid the construction of the strict super-solutions that were required in TB22. Another difference is that the nonlocal terms are nonlinear in the former, while are linear in the latter. Nevertheless, we can use their methodology because the jumps are compound Poisson types having bounded variations.

It is sufficient to show the inequality

$$\bar{F} \geq \underline{F} \quad \text{in} \quad \bar{\Omega} \tag{80}$$

for any super-solutions $\bar{F}$ and sub-solutions $\underline{F}$. Suppose that (80) does not hold true. Then, there exists a point $(x_1^*, x_2^*) \in \bar{\Omega}$ such that

$$M := \underline{F}(x_1^*, x_2^*) - \bar{F}(x_1^*, x_2^*) > 0. \tag{81}$$

If $(x_1^*, x_2^*) \in \Omega$, then the proof follows a classical argument based on the standard technique of the doubling variables (e.g., use the auxiliary function of Eq. (5.15) without the last term in [71]). We therefore deal with the boundary case $(x_1^*, x_2^*) \in \partial\Omega$.

Because $\Omega$ is a square and is therefore a polygonal domain, there exist uniformly continuous functions $\gamma_1, \gamma_2 : \bar{\Omega} \to \mathbb{R}$ such that

$$B\big((x_1, x_2) + h(\gamma_1(x_1, x_2), \gamma_1(x_1, x_2)), hr\big) \subset \Omega \quad \text{and} \quad h \in (0, h_0] \tag{82}$$

where $B(X, R)$ is the 2D open ball with the center $X$ and radius $R > 0$, and $h_0, r > 0$ are sufficiently small constant. For any $\kappa > 1$ and $\mu \in (0,1)$, we set the auxiliary function $\Theta : \bar{\Omega} \times \bar{\Omega} \to \mathbb{R}$:

$$\begin{aligned}\Theta(x_1, x_2, y_1, y_2) = &\underline{F}(x_1, x_2) - \bar{F}(y_1, y_2) \\ &- \left\{ \left(\kappa(x_1 - y_1) + \mu\gamma_1(x_1^*, x_2^*)\right)^2 + \left(\kappa(x_2 - y_2) + \mu\gamma_2(x_1^*, x_2^*)\right)^2 \right\} \\ &- \mu\left((x_1 - x_1^*)^2 + (x_2 - x_2^*)^2\right)\end{aligned} \tag{83}$$

We set $M_\kappa = \arg\max\limits_{(x_1, x_2, y_1, y_2) \in \bar{\Omega} \times \bar{\Omega}} \Theta(x_1, x_2, y_1, y_2)$. Following the proof of Theorem 5.1 of TB22, the assumption $\underline{F}, \bar{F} \in BUC(\bar{\Omega})$ leads to that a maximizer $(\hat{x}_1, \hat{x}_2, \hat{y}_1, \hat{y}_2) \in \bar{\Omega} \times \bar{\Omega}$ of $\Theta$, which depends on $\kappa$ and $\mu$, exists and satisfies the following properties by choosing a suitable subsequence that will be assumed herein:

$$\left|\kappa(\hat{x}_1 - \hat{y}_1)\right|, \left|\kappa(\hat{x}_2 - \hat{y}_2)\right| \quad \text{are uniformly bounded for} \quad \kappa > 1 \tag{84}$$

and hence

$$\lim_{\kappa \to +\infty}(\hat{x}_1 - \hat{y}_1) = \lim_{\kappa \to +\infty}(\hat{x}_2 - \hat{y}_2) = 0, \tag{85}$$

and further

$$\lim_{\kappa \to +\infty} \left( \underline{F}(\hat{x}_1, \hat{x}_2) - \overline{F}(\hat{y}_1, \hat{y}_2) \right) \leq M \tag{86}$$

as well as

$$\lim_{\kappa \to +\infty} \left( \kappa(\hat{x}_1 - \hat{y}_1) + \mu \gamma_1(x_1^*, x_2^*) \right) = \lim_{\kappa \to +\infty} \left( \kappa(\hat{x}_2 - \hat{y}_2) + \mu \gamma_2(x_1^*, x_2^*) \right) = 0. \tag{87}$$

Consequently, we arrive at

$$\lim_{\kappa \to +\infty} (\hat{x}_1, \hat{y}_1) = \lim_{\kappa \to +\infty} (\hat{x}_2, \hat{y}_2) = (x_1^*, x_2^*) \text{ and } \lim_{\kappa \to +\infty} M_\kappa = M. \tag{88}$$

By utilizing the uniform continuity of the functions $\gamma_1, \gamma_2$, we obtain

$$\hat{y}_1 = \hat{x}_1 + \frac{\mu}{\kappa} \gamma_1(\hat{x}_1, \hat{x}_2) + o(\kappa^{-1}) \text{ and } \hat{y}_2 = \hat{x}_2 + \frac{\mu}{\kappa} \gamma_2(\hat{x}_1, \hat{x}_2) + o(\kappa^{-1}) \tag{89}$$

with $o(\cdot)$ the classical Landau symbol. The relationship (89) guarantees $(\hat{y}_1, \hat{y}_2) \in \Omega$ for sufficiently large $\kappa$. Then, it follows that the following $\overline{\varphi}, \underline{\varphi} \in C^1(\overline{\Omega})$ satisfy that $\overline{F} - \overline{\varphi}$ is maximized at $(\hat{x}_1, \hat{x}_2)$ and $\underline{F} - \underline{\varphi}$ is minimized at $(\hat{y}_1, \hat{y}_2)$:

$$\begin{aligned}\overline{\varphi}(x_1, x_2) &= \overline{F}(\hat{y}_1, \hat{y}_2) \\ &+ \left( \kappa(x_1 - \hat{y}_1) + \mu \gamma_1(x_1^*, x_2^*) \right)^2 + \left( \kappa(x_2 - \hat{y}_2) + \mu \gamma_2(x_1^*, x_2^*) \right)^2 \\ &+ \mu \left( (x_1 - x_1^*)^2 + (x_2 - x_2^*)^2 \right) \end{aligned} \tag{90}$$

$$\begin{aligned}\underline{\varphi}(y_1, y_2) &= \underline{F}(\hat{x}_1, \hat{x}_2) \\ &- \left\{ \left( \kappa(\hat{x}_1 - y_1) + \mu \gamma_1(x_1^*, x_2^*) \right)^2 + \left( \kappa(\hat{x}_2 - y_2) + \mu \gamma_2(x_1^*, x_2^*) \right)^2 \right\}. \\ &- \mu \left( (\hat{x}_1 - x_1^*)^2 + (\hat{x}_2 - x_2^*)^2 \right) \end{aligned} \tag{91}$$

By Definition 1, we obtain the inequalities

$$\begin{aligned}&\lambda_Z \inf_{(u_1, u_2) \in C \times C} \left\{ \overline{F}(\hat{x}_1 + c_1(\hat{x}_1, \hat{x}_2, u_1), \hat{x}_2 + c_2(\hat{x}_1, \hat{x}_2, u_2)) - \overline{F}(\hat{x}_1, \hat{x}_2) + g^p(u_1, u_2) \right\} \\ &- \left( \delta + \hat{\delta}(\overline{F})(\hat{x}_1, \hat{x}_2) \right) \overline{F}(\hat{x}_1, \hat{x}_2) + a_1(\hat{x}_1, \hat{x}_2) \frac{\partial \overline{\varphi}(\hat{x}_1, \hat{x}_2)}{\partial x_1} + a_2(\hat{x}_1, \hat{x}_2) \frac{\partial \overline{\varphi}(\hat{x}_1, \hat{x}_2)}{\partial x_2} + f^p(\hat{x}_1, \hat{x}_2) \leq 0 \end{aligned} \tag{92}$$

and

$$\begin{aligned}&\lambda_Z \inf_{(u_1, u_2) \in C \times C} \left\{ \underline{F}(\hat{x}_1 + c_1(\hat{x}_1, \hat{x}_2, u_1), \hat{x}_2 + c_2(\hat{x}_1, \hat{x}_2, u_2)) - \underline{F}(\hat{x}_1, \hat{x}_2) + g^p(u_1, u_2) \right\} \\ &- \left( \delta + \hat{\delta}(\underline{F})(\hat{x}_1, \hat{x}_2) \right) \underline{F}(\hat{x}_1, \hat{x}_2) + a_1(\hat{x}_1, \hat{x}_2) \frac{\partial \underline{\varphi}(\hat{x}_1, \hat{x}_2)}{\partial x_1} + a_2(\hat{x}_1, \hat{x}_2) \frac{\partial \underline{\varphi}(\hat{x}_1, \hat{x}_2)}{\partial x_2} + f^p(\hat{x}_1, \hat{x}_2) \geq 0 \end{aligned}. \tag{93}$$

Summing up (92) and $-1 \times$ (93), and then take limits of $\kappa \to +\infty$ and $\mu \to +0$ in this order, to obtain

$$\begin{aligned}&-\lambda_Z \inf_{(u_1, u_2) \in C \times C} \left\{ \underline{F}(x_1^* + c_1(x_1^*, x_2^*, u_1), x_2^* + c_2(x_1^*, x_2^*, u_2)) - \underline{F}(x_1^*, x_2^*) + g^p(u_1, u_2) \right\} \\ &+ \lambda_Z \inf_{(u_1, u_2) \in C \times C} \left\{ \overline{F}(x_1^* + c_1(x_1^*, x_2^*, u_1), x_2^* + c_2(x_1^*, x_2^*, u_2)) - \overline{F}(x_1^*, x_2^*) + g^p(u_1, u_2) \right\} \\ &+ \left( \hat{\delta}(\underline{F})(x_1^*, x_2^*) \right) \underline{F}(x_1^*, x_2^*) - \left( \hat{\delta}(\overline{F})(x_1^*, x_2^*) \right) \overline{F}(x_1^*, x_2^*) \\ &+ \delta \left( \underline{F}(x_1^*, x_2^*) - \overline{F}(x_1^*, x_2^*) \right) \leq 0 \end{aligned}. \tag{94}$$

We have the inequality ( $(\underline{u}_1, \underline{u}_2)$ is a minimizer of the first line that exists because $\underline{F} \in BUC(\bar{\Omega})$, $C$ is compact, and $g$ is lower-semicontinuous in $C \times C$ )

$$
\begin{aligned}
&-\inf_{(u_1,u_2)\in C\times C}\left\{\underline{F}\left(x_1^* + c_1(x_1^*, x_2^*, u_1), x_2^* + c_2(x_1^*, x_2^*, u_2)\right) - \underline{F}(x_1^*, x_2^*) + g^p(u_1, u_2)\right\} \\
&+\inf_{(u_1,u_2)\in C\times C}\left\{\bar{F}\left(x_1^* + c_1(x_1^*, x_2^*, u_1), x_2^* + c_2(x_1^*, x_2^*, u_2)\right) - \bar{F}(x_1^*, x_2^*) + g^p(u_1, u_2)\right\} \\
&= -\min_{(u_1,u_2)\in C\times C}\left\{\underline{F}\left(x_1^* + c_1(x_1^*, x_2^*, u_1), x_2^* + c_2(x_1^*, x_2^*, u_2)\right) - \underline{F}(x_1^*, x_2^*) + g^p(u_1, u_2)\right\} \\
&+\min_{(u_1,u_2)\in C\times C}\left\{\bar{F}\left(x_1^* + c_1(x_1^*, x_2^*, u_1), x_2^* + c_2(x_1^*, x_2^*, u_2)\right) - \bar{F}(x_1^*, x_2^*) + g^p(u_1, u_2)\right\} \\
&\geq -\underline{F}\left(x_1^* + c_1(x_1^*, x_2^*, \underline{u}_1), x_2^* + c_2(x_1^*, x_2^*, \underline{u}_2)\right) + \underline{F}(x_1^*, x_2^*) - g^p(\underline{u}_1, \underline{u}_2) \\
&+\bar{F}\left(x_1^* + c_1(x_1^*, x_2^*, \underline{u}_1), x_2^* + c_2(x_1^*, x_2^*, \underline{u}_2)\right) - \bar{F}(x_1^*, x_2^*) + g^p(\underline{u}_1, \underline{u}_2) \\
&= \underline{F}(x_1^*, x_2^*) - \bar{F}(x_1^*, x_2^*) \\
&\quad -\left\{\underline{F}\left(x_1^* + c_1(x_1^*, x_2^*, \underline{u}_1), x_2^* + c_2(x_1^*, x_2^*, \underline{u}_2)\right) - \bar{F}\left(x_1^* + c_1(x_1^*, x_2^*, \underline{u}_1), x_2^* + c_2(x_1^*, x_2^*, \underline{u}_2)\right)\right\} \\
&\geq \underline{F}(x_1^*, x_2^*) - \bar{F}(x_1^*, x_2^*) - \left\{\underline{F}(x_1^*, x_2^*) - \bar{F}(x_1^*, x_2^*)\right\} \\
&= 0
\end{aligned}
\qquad (95)
$$

The last line is due to the maximizing property of $(x_1^*, x_2^*)$. We have

$$
\begin{aligned}
&\frac{\psi}{\lambda_N} \times \left\{\left(\hat{\delta}(\underline{F})(x_1^*, x_2^*)\right)\underline{F}(x_1^*, x_2^*) - \left(\hat{\delta}(\bar{F})(x_1^*, x_2^*)\right)\bar{F}(x_1^*, x_2^*)\right\} \\
&= \int_{-1}^{+\infty}\left[1 - \exp\left(\psi\omega(\underline{F})(x_1^*, x_2^*, z)\right)\right] P(z)\,dz\,\underline{F}(x_1^*, x_2^*) \\
&\quad - \int_{-1}^{+\infty}\left[1 - \exp\left(\psi\omega(\bar{F})(x_1^*, x_2^*, z)\right)\right] P(z)\,dz\,\bar{F}(x_1^*, x_2^*) \\
&= \int_{-1}^{+\infty}\left[1 - \exp\left(\psi\frac{\underline{F}\left(x_1^* + b_1(x_1^*, x_2^*)z, x_2^* + b_2(x_1^*, x_2^*)z\right) - \underline{F}(x_1^*, x_2^*)}{\underline{F}(x_1^*, x_2^*)}\right)\right] P(z)\,dz\,\underline{F}(x_1^*, x_2^*) \\
&\quad - \int_{-1}^{+\infty}\left[1 - \exp\left(\psi\frac{\bar{F}\left(x_1^* + b_1(x_1^*, x_2^*)z, x_2^* + b_2(x_1^*, x_2^*)z\right) - \bar{F}(x_1^*, x_2^*)}{\bar{F}(x_1^*, x_2^*)}\right)\right] P(z)\,dz\,\bar{F}(x_1^*, x_2^*)
\end{aligned}
\qquad (96)
$$

Due again to the maximizing property of $(x_1^*, x_2^*)$, we have

$$
\begin{aligned}
\underline{F}(x_1^*, x_2^*) - \bar{F}(x_1^*, x_2^*) &\geq \underline{F}\left(x_1^* + b_1(x_1^*, x_2^*)z, x_2^* + b_2(x_1^*, x_2^*)z\right) \\
&\quad - \bar{F}\left(x_1^* + b_1(x_1^*, x_2^*)z, x_2^* + b_2(x_1^*, x_2^*)z\right)
\end{aligned}
\qquad (97)
$$

inside the integrations of (96), and hence

$$
\begin{aligned}
&\bar{F}\left(x_1^* + b_1(x_1^*, x_2^*)z, x_2^* + b_2(x_1^*, x_2^*)z\right) - \bar{F}(x_1^*, x_2^*) \\
&\geq \underline{F}\left(x_1^* + b_1(x_1^*, x_2^*)z, x_2^* + b_2(x_1^*, x_2^*)z\right) - \underline{F}(x_1^*, x_2^*)
\end{aligned},
\qquad (98)
$$

to obtain (See, Lemma 1 for $y$)

$$\int_{-1}^{+\infty}\left[1-\exp\left(\psi\frac{\underline{F}\left(x_1^*+b_1\left(x_1^*,x_2^*\right)z,x_2^*+b_2\left(x_1^*,x_2^*\right)z\right)-\underline{F}\left(x_1^*,x_2^*\right)}{\underline{F}\left(x_1^*,x_2^*\right)}\right)\right]P(z)\mathrm{d}z\underline{F}\left(x_1^*,x_2^*\right)$$
$$-\int_{-1}^{+\infty}\left[1-\exp\left(\psi\frac{\overline{F}\left(x_1^*+b_1\left(x_1^*,x_2^*\right)z,x_2^*+b_2\left(x_1^*,x_2^*\right)z\right)-\overline{F}\left(x_1^*,x_2^*\right)}{\overline{F}\left(x_1^*,x_2^*\right)}\right)\right]P(z)\mathrm{d}z\overline{F}\left(x_1^*,x_2^*\right)$$
$$\geq \int_{-1}^{+\infty}\left[1-\exp\left(\psi\frac{\underline{F}\left(x_1^*+b_1\left(x_1^*,x_2^*\right)z,x_2^*+b_2\left(x_1^*,x_2^*\right)z\right)-\underline{F}\left(x_1^*,x_2^*\right)}{\underline{F}\left(x_1^*,x_2^*\right)}\right)\right]P(z)\mathrm{d}z\underline{F}\left(x_1^*,x_2^*\right) \quad (99)$$
$$-\int_{-1}^{+\infty}\left[1-\exp\left(\psi\frac{\underline{F}\left(x_1^*+b_1\left(x_1^*,x_2^*\right)z,x_2^*+b_2\left(x_1^*,x_2^*\right)z\right)-\underline{F}\left(x_1^*,x_2^*\right)}{\overline{F}\left(x_1^*,x_2^*\right)}\right)\right]P(z)\mathrm{d}z\overline{F}\left(x_1^*,x_2^*\right)$$

Due to $\overline{F}(x_1^*,x_2^*) < \underline{F}(x_1^*,x_2^*)$, we further have (See, Lemma 1 for $x$)

$$\int_{-1}^{+\infty}\left[1-\exp\left(\psi\frac{\underline{F}\left(x_1^*+b_1\left(x_1^*,x_2^*\right)z,x_2^*+b_2\left(x_1^*,x_2^*\right)z\right)-\underline{F}\left(x_1^*,x_2^*\right)}{\underline{F}\left(x_1^*,x_2^*\right)}\right)\right]P(z)\mathrm{d}z\underline{F}\left(x_1^*,x_2^*\right)$$
$$-\int_{-1}^{+\infty}\left[1-\exp\left(\psi\frac{\underline{F}\left(x_1^*+b_1\left(x_1^*,x_2^*\right)z,x_2^*+b_2\left(x_1^*,x_2^*\right)z\right)-\underline{F}\left(x_1^*,x_2^*\right)}{\overline{F}\left(x_1^*,x_2^*\right)}\right)\right]P(z)\mathrm{d}z\overline{F}\left(x_1^*,x_2^*\right)$$
$$\geq \int_{-1}^{+\infty}\left[1-\exp\left(\psi\frac{\underline{F}\left(x_1^*+b_1\left(x_1^*,x_2^*\right)z,x_2^*+b_2\left(x_1^*,x_2^*\right)z\right)-\underline{F}\left(x_1^*,x_2^*\right)}{\overline{F}\left(x_1^*,x_2^*\right)}\right)\right]P(z)\mathrm{d}z\overline{F}\left(x_1^*,x_2^*\right). \quad (100)$$
$$-\int_{-1}^{+\infty}\left[1-\exp\left(\psi\frac{\underline{F}\left(x_1^*+b_1\left(x_1^*,x_2^*\right)z,x_2^*+b_2\left(x_1^*,x_2^*\right)z\right)-\underline{F}\left(x_1^*,x_2^*\right)}{\overline{F}\left(x_1^*,x_2^*\right)}\right)\right]P(z)\mathrm{d}z\overline{F}\left(x_1^*,x_2^*\right)$$
$$\geq 0$$

By (95) and (99)-(100), (94) leads to

$$\delta\left(\underline{F}\left(x_1^*,x_2^*\right)-\overline{F}\left(x_1^*,x_2^*\right)\right)=\delta M\leq 0, \quad (101)$$

which is a contradiction. Therefore, $M$ in (81) must not be positive, and hence (80) holds true.

□

**Lemma 1** *Define the function* $\rho:(0,+\infty)\times\mathbb{R}\to\mathbb{R}$ *as*

$$\rho(x,y)=x\left(1-\exp\left(\frac{y}{x}\right)\right). \quad (102)$$

*This* $\rho$ *is increasing and decreasing with respect to* $x\in(0,+\infty)$ *and* $y\in\mathbb{R}$, *respectively.*

**Proof of Lemma 1**

The proof is by the direct partial differentiations as follows:

$$\frac{\partial}{\partial x}\rho(x,y)=1-\exp\left(\frac{y}{x}\right)+x\frac{y}{x^2}\exp\left(\frac{y}{x}\right)=\left(\exp\left(-\frac{y}{x}\right)-\left(-\frac{y}{x}\right)-1\right)\exp\left(\frac{y}{x}\right)\geq 0 \quad (103)$$

due to the elementary inequality

$$\exp(x) - x - 1 \geq 0, \quad x \in \mathbb{R}, \tag{104}$$

and

$$\frac{\partial}{\partial y} \rho(x, y) = -\exp\left(\frac{y}{x}\right) < 0. \tag{105}$$

□

*Proof of Proposition 4*

Consider firstly the auxiliary system

$$G'_{i,j}(F) = f^p(ih_1, jh_2), \quad i = 0, 1, 2, \ldots, N_1, \quad j = 0, 1, 2, \ldots, N_2, \tag{106}$$

where (for colored terms, see **Section 3** in the main text)

$$\begin{aligned}
G'_{i,j}(F) := &\, \overline{C}_{i,j} \left( \delta + \lambda_Z + \frac{\lambda_N}{\psi} \sum_{m=1}^{M} \exp\left( -\psi \frac{\min\{\omega, \max\{-\omega, F_{i,j} - \hat{F}_{i,j,m}\}\}}{\max\{F_{i,j}, \varepsilon\}} \right) p_m \right) \max\{F_{i,j}, 0\} \\
&+ \left( \begin{array}{l} \dfrac{a_1(ih_1, jh_2)}{h_1} \times \begin{cases} F_{i,j} - F_{i+1,j} & (a_1(ih_1, jh_2) \geq 0) \\ F_{i-1,j} - F_{i,j} & (a_1(ih_1, jh_2) < 0) \end{cases} \\ + \dfrac{a_2(ih_1, jh_2)}{h_2} \times \begin{cases} F_{i,j} - F_{i,j+1} & (a_2(ih_1, jh_2) \geq 0) \\ F_{i,j-1} - F_{i,j} & (a_2(ih_1, jh_2) < 0) \end{cases} \end{array} \right) \\
&- \lambda_Z \min_{(k,l) \in \{0,1,2,\ldots,N_1\} \times \{0,1,2,\ldots,N_2\}} \left\{ F_{i-\min\{i,k\}, j-\min\{j,l\}} - F_{i,j} + g^p(kh_1, lh_2) \right\}
\end{aligned} \tag{107}$$

with a constant $\omega > 0$ chosen later. The difference between $G$ and $G'$ is that the first lines of (64) and (107) such that the latter uses "max" and the truncation using $\omega > 0$. With this truncation, it follows that $G'_{i,j}(F)$ is increasing with respect to the green "$F_{i,j}$"s, and further it is also increasing with respect to the magenta "$F_{i,j} - F_{\cdot,\cdot}$"s irrespective to the sign of $F_{i,j}$. Moreover, it follows that (106) satisfies the Lipschitz continuity and is proper due to $\delta > 0$ (Definitions 3 and 6 of Oberman [78]). We then directly apply Theorem 8 of Oberman [78] to obtain that the auxiliary system (106) admits a unique numerical solution (63), denoted as $\hat{F} = \{\hat{F}_{i,j}\}$, such that (65) is satisfied. Notice that the numerical solution $\hat{F}$ depends on $\omega$, while it also satisfies (63) if $\omega > 0$ is sufficiently large such that

$$\omega > \underline{\omega} := \frac{1}{\delta} \left( \max_{(i,j) \in \{0,1,2,\ldots,N_1\} \times \{0,1,2,\ldots,N_2\}} f^p(ih_1, jh_2) - \min_{(i,j) \in \{0,1,2,\ldots,N_1\} \times \{0,1,2,\ldots,N_2\}} f^p(ih_1, jh_2) \right). \tag{108}$$

Further, the last "max" in the first line of (107) is redundant if $F = \hat{F}$. Hereafter, we assume $\omega > \underline{\omega}$.

We finally show that $F = \hat{F}$ is the unique solution to (63) such that (65). If it is not true, then there exists more than one numerical solutions to (63). Then, we arrive at the conclusion that the auxiliary system (107) admits more than one numerical solutions complying with (65), but this is false due to the result obtained in the previous paragraph. Consequently, there exists only one numerical solution to (63) such that (65), and the proposition is proven.

**Appendix B: Field survey results**

Field survey results not presented in Section 4 are provided in this appendix. **Table B1** shows the observed populations evaluated in terms of the area ($m^2$) in the observation strip at each Station. **Table B2** then shows the normalized populations obtained by dividing the area in **Table B1** by the total area 480 ($m^2$), which were used in the model application. **Table B3** shows the best-fit results of the normalized population by the proposed model. **Table B4** provides the relative errors between the observed and fitted results. The relative error averaged against all the observation data (excluding April 21 at which the initial condition coincides with the observation value) is at most 12.4% at Station 3 and is the smallest value 0.6% at Station 7.

**Table B1.** Observed populations ($m^2$) at each Station in 2022.

| Date in 2022 | April 21 | May 19 | June 16 | July 28 | August 23 | September 15 |
|---|---|---|---|---|---|---|
| Station 1 | 2 | 5 | 122 | 19 | 52 | 1 |
| Station 2 | 0 | 1 | 37 | 47 | 79 | 24 |
| Station 3 | 0 | 0 | 41 | 107 | 367 | 18 |
| Station 4 | 0 | 0 | 8 | 91 | 13 | 105 |
| Station 5 | 0 | 0 | 1 | 65 | 119 | 69 |
| Station 6 | 0 | 0 | 8 | 151 | 183 | 13 |
| Station 7 | 0 | 3 | 6 | 121 | 335 | 3 |

**Table B2.** Normalized observed populations at each Station in 2022.

| Date in 2022 | April 21 | May 19 | June 16 | July 28 | August 23 | September 15 |
|---|---|---|---|---|---|---|
| Station 1 | 0.0042 | 0.0104 | 0.2542 | 0.0396 | 0.1083 | 0.0021 |
| Station 2 | 0.0000 | 0.0021 | 0.0771 | 0.0979 | 0.1646 | 0.0500 |
| Station 3 | 0.0000 | 0.0000 | 0.0854 | 0.2229 | 0.7646 | 0.0375 |
| Station 4 | 0.0000 | 0.0000 | 0.0167 | 0.1896 | 0.0271 | 0.2188 |
| Station 5 | 0.0000 | 0.0000 | 0.0021 | 0.1354 | 0.2479 | 0.1438 |
| Station 6 | 0.0000 | 0.0000 | 0.0167 | 0.3146 | 0.3813 | 0.0271 |
| Station 7 | 0.0000 | 0.0063 | 0.0125 | 0.2521 | 0.6979 | 0.0063 |

**Table B3.** Fitted populations at each Station in 2022.

| Date in 2022 | April 21 | May 19 | June 16 | July 28 | August 23 | September 15 |
|---|---|---|---|---|---|---|
| Station 1 | 0.0042 | 0.0162 | 0.0296 | 0.0517 | 0.0675 | 0.0827 |
| Station 2 | 0.0000 | 0.0410 | 0.0684 | 0.0943 | 0.1009 | 0.1044 |
| Station 3 | 0.0000 | 0.2453 | 0.2546 | 0.1420 | 0.6764 | 0.0000 |
| Station 4 | 0.0000 | 0.0140 | 0.0349 | 0.0816 | 0.1299 | 0.1883 |
| Station 5 | 0.0000 | 0.0257 | 0.0613 | 0.1278 | 0.1701 | 0.1999 |
| Station 6 | 0.0000 | 0.0331 | 0.0971 | 0.2343 | 0.2409 | 0.2393 |
| Station 7 | 0.0000 | 0.0048 | 0.0322 | 0.2504 | 0.7005 | 0.0000 |

**Table B4.** Relative errors between the observed and fitted results.

| Date in 2022 | April 21 | May 19 | June 16 | July 28 | August 23 | September 15 |
|---|---|---|---|---|---|---|
| Station 1 | 0.000 | 0.006 | 0.225 | 0.012 | 0.041 | 0.081 |
| Station 2 | 0.000 | 0.039 | 0.009 | 0.004 | 0.064 | 0.054 |
| Station 3 | 0.000 | 0.245 | 0.169 | 0.081 | 0.088 | 0.038 |
| Station 4 | 0.000 | 0.014 | 0.018 | 0.108 | 0.103 | 0.030 |
| Station 5 | 0.000 | 0.026 | 0.059 | 0.008 | 0.078 | 0.056 |
| Station 6 | 0.000 | 0.033 | 0.080 | 0.080 | 0.140 | 0.212 |
| Station 7 | 0.000 | 0.001 | 0.020 | 0.002 | 0.003 | 0.006 |

**Appendix C: On the computational resolution**

The computational resolution used in our numerical analysis is sufficiently fine to obtain the conclusion of Section 4, as demonstrated in **Figure C.1** below where the computed optimal $u_1$ at each Station 2 against different resolutions are presented. Although not presented here, the value function and the other control variable are comparable among these resolutions as well, supporting the use of the resolution.

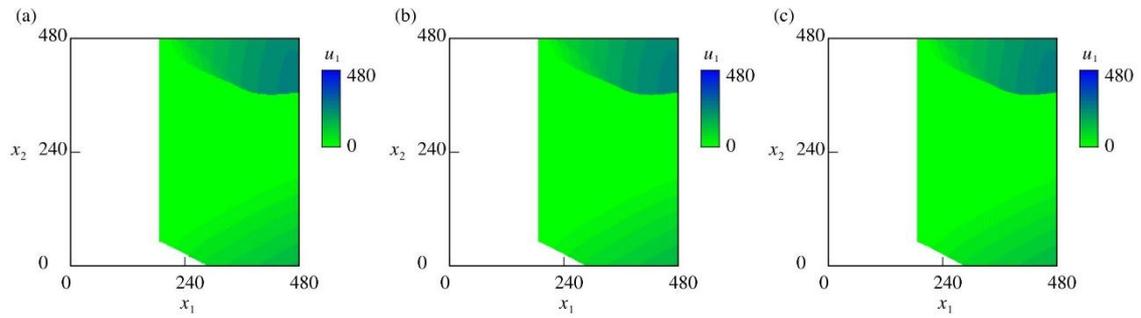

**Figure C.1.** The computed optimal $u_1$ at each Station 2 where the domain is uniformly discretized with (a) 201, (b) 301, and (c) 401 vertices in each direction. The white area represents $u_1 = 0$ (Do nothing against the drifting population).

**Competing interests** The authors declare no competing interests.

**Acknowledgements** None

**Funding** This work was supported by the Japan Society for the Promotion of Science (No. 22K14441) and the Ministry of Land, Infrastructure, and Transport of Japan (B4R202003, B4R202101, and B4R202201).

**Availability of data and material** Data are available upon reasonable request to the corresponding author.